\documentclass[14pt]{article}
\usepackage[english]{babel}
\usepackage{todonotes}
\usepackage{amsmath,amsfonts,amsthm,bm,dsfont} 
\usepackage[margin=3cm]{geometry}
\usepackage{enumerate}
\usepackage{setspace}
\newtheorem{theorem}{Theorem}

\newtheorem{corollary}[theorem]{Corollary}

\newtheorem{example}[theorem]{Example}

\newtheorem{lemma}[theorem]{Lemma}

\newtheorem{proposition}[theorem]{Proposition}

\numberwithin{theorem}{section}
\newcommand{\Rd}{\mathbb{R}^d}

\begin{document}

\title{The identification of diffusions from imperfect observations}

\author{J.M.C.Clark \and Dan Crisan}
\date{22-07-2024}
\maketitle

\begin{abstract}

This paper studies the identification of an $\Rd$-valued diffusion $X$ when a running function of it, say $h(X_t)$, is observed. A point-wise observation of the process (in other words, observing $h(X_t)$ in isolation) cannot identify $X_t$ unless the $h$ is injective. However observing $h(X_s)$ on a small interval $[t,t+\varepsilon]$ can be enough to determine $X_t$ exactly.
The paper contain results that expand on this idea; in particular, a property of `fine total asymmetry' of twice continuously differentiable $h$ is introduced that depends on the fine topology of potential theory and that is both necessary and sufficient for $X$ to be adapted to a natural right-continuous filtration generated by the observations. This particular filtration, though augmented with null sets, does not depend on the distribution of $X_0$.  For real-analytic $h$ the property reduces to simple asymmetry; that is, there is no nontrivial affine isometry $\kappa$ on $\Rd$ such that $h = h \circ \kappa$.
A second result concerns the case where $X_0$ is given and $h$ is merely Borel;  then  
 $X$ is adapted to an augmented filtration generated by the observation process $(h(X_t))_{t\geq 0}$ if $h$ is `locally invertible' on a subset of $\Rd$ dense in the fine topology on $\Rd$. 

\end{abstract}

\section{Introduction}

Consider the following question about the measurability of a process. Suppose a continuous stochastic process is continuously measured through a nonlinearity that is not one-to-one. Under what conditions is the process `immediately and continuously determined' in the sense that it is adapted to the right-continuous filtration generated, with suitable completion, by the process of measurements? For the simple case of a real-valued Brownian motion with a strictly positive initial density that is measured through a real-analytic function h, the answer turns out to be straightforward: the Brownian motion is adapted to the measured filtration if and only if $h$ is aperiodic and $h(a + \cdot)$ is not an even function for any real value $a$, or, in other words, if and only if $h$ is asymmetric in the sense that it is not invariant under any non-trivial Euclidean isometry or ``rigid motion'' on the real line. Other cases are not as clear-cut. The purpose of this paper is to present this and more general results of this type for a class of multidimensional processes similar to Brownian motion and for various classes of measurement maps. 
 
A broader version of the question just posed splits naturally into two separate problems, which might loosely be termed the `global' problem of immediate determination and the `local' problem of continuous `tracking'. 
Suppose $X = (X_t)_{t \geq 0}$ is an $\Rd$-valued Brownian motion with an initial distribution $\mu$ and $h: \Rd \rightarrow \mathbb{R}^e$ is the Borel measurement map defining an observation process $Y = (h(X_t))_{t \geq 0}$. This paper presents conditions on $\mu$ and $h$ under which: 

\noindent (a) $X_0$ is determined by -- that is, measurable with respect to --  the initial ``germ field'' ${\cal Y}^\mu_{0+}$ of $Y$; 

\noindent (b) the process $X$ is adapted to the augmented right-continuous filtration $({\cal Y}^\mu_{t+})_{t \geq 0}$\footnote{Here, for $ t \geq 0, \;{\cal Y}^\mu_{t+} = \bigcap_{t \leq s}\sigma(Y_r:\:0\leq r \leq s; \:{\cal N}^\mu_s),\;{\cal N}^\mu_s$ being a family of null sets. More precise definitions are given in the next section.} generated by $Y$.
 
 The results extend to other continuous processes that are closely related to Brownian motion, for instance in the sense that they generate the same right-continuous locally completed filtration on the space of $\Rd$-valued continuous paths, or are transformed into Brownian motion by a homeomorphic change of coordinates. The fine topology of potential theory plays a central role in the formulation of the conditions.
The main result is that, for a $C^2(\Rd,\Rd)$ measurement map $h$ with a Jacobian that is  non-singular on an open finely dense set and a distribution of $X_0$ that charges all finely open sets, $X_0$ is determined by the initial germ field of the observed process $Y$ if and only if $h$ is ``finely totally asymmetric'', a notion of local asymmetry that depends on the fine topology.
A further result is that $X$ is adapted to an augmented filtration generated by $Y$ if $X_0$ is given and $h$ is `locally invertible' on a finely dense subset of $\Rd$

\indent The problems posed in the paragraphs above have their origins in certain singular formulations of the non-linear filtering problem of diffusion processes, in which the conditional distribution $\hat{\mu}_t = \mathds{P}(X_t \in \cdot |{\cal Y}_{t+} ) $ of the current value of a ``signal'' diffusion $X$, conditioned on the past history of a related observation process $Y$, 
turns out to be a random Dirac measure $\delta_{X_t}$. This signifies that  $X_t$ is observable in the sense that it is determined exactly by the past of $Y$ with probability one. Such examples occur in the ``signal dependent observation noise'' models considered in \cite{joannides1997nonlinear} and \cite{crisan2009nonlinear}. One simplified version illustrates the point. Consider a real-valued observation process  $Y$ that is related to an $\mathbb{R}^d$-valued diffusion $X$ via an It\^o equation of the form
$dY_t = g(X_t)dV_t$,
where  $V$ is a real-valued Brownian motion independent of $X$ and $g$ is a continuous function. It follows from standard theorems on quadratic variation (see, e.g.  \cite{protter2003stochastic}, Chap. II) that the quadratic variation process $ \int_0^t g(X_s)^2ds$ of $Y$, along with the continuous derivative process $g(X_t)^2$, are both adapted to the observed right-continuous filtration generated (with suitable augmentation) by $Y$. A natural question to ask is under what conditions on $g$ is $X$ similarly adapted to this observed filtration.

Separately, in a version of the ``small additive noise'' problem, an $\mathbb{R}^e$-valued process $Y^\epsilon$ of observations of a $\Rd$-valued diffusion $X$ is modelled by an It\^o equation  containing a small parameter $\varepsilon$:
$ dY_t^\epsilon  = h(X_t)dt + \varepsilon dV_t$. Here $V$ is a $\Rd$-valued Brownian motion, independent of $X$, representing low-intensity noise. In some formulations, $X$ also depends on $\varepsilon$ and is modelled by an It\^o differential equation that degenerates to an ordinary one as $\varepsilon$ vanishes.        
This problem, important in applications, has had an asymptotic treatment in a number of papers, dating back to \cite{katzur1984asymptotic1},  \cite{krener1986asymptotic} and \cite{picard1986nonlinear}, where it is shown that various approximate filtering schemes produce  estimates $\hat{X}_t^\varepsilon$ of $X_t$ that are adapted to the observable filtration generated by $Y^\epsilon$ and that converge to $X_t$, at least in probability, as $\varepsilon$ vanishes. In most of these papers $h$ is assumed to be one-to-one; relatively few consider cases where $h$ is many-to-one. \cite{krener1986asymptotic} contains interesting results on multiple time-scales for the case where $e<d$, $h$ is linear and $X$ degenerates with $\varepsilon$. The papers most relevant to this paper are \cite{fleming1988piecewise} and \cite{fleming1989piecewise}, which treat cases where $d=e=1$, the diffusion is independent of $\varepsilon$ and $h$ is, respectively, a piecewise linear and a piecewise monotone many-to-one function. In the last two papers the authors essentially solve the global problem for the limiting form (with no `additive noise') of the cases they consider. Generalisations of this `noiseless' idealization of their results are presented in the final section.

\section{Preliminaries}
In what follows $(\Omega,{\cal X}^0)$ is a canonical measurable space for an $\mathbb{R}^d$-valued continuous process $X$: that is, $\Omega=C([0,\infty),\mathbb{R}^d)$ and ${\cal X}^0$ is the $\sigma$-algebra generated by the coordinate maps $X_t(\omega) = \omega(t)$. Let $\mathbb{X}^0$ denote $({\cal X}^0_{t+})_{t \geq 0}$, the natural right-continuous filtration of $X$, where ${\cal X}^0_{t} = \sigma(X_r: 0 \leq r \leq t) \text{ and }  {\cal X}^0_{t+} = \bigcap\,_{s>t} {\cal X}^0_{s}  $. Unless stated otherwise we assume that the process $X = (X_t)_{t\geq 0}$ is a diffusion process that is `locally equivalent in law' to Brownian motion in the following sense: $(\Omega,{\cal X}^0)$ is endowed with a strong Markov family of laws $\{\widetilde{\mathds{P}}^x\}_{x \in \Rd}$, each corresponding to the law of $X$ starting at $x$, such that, for every $t$, each $\widetilde{\mathds{P}}^x$ is equivalent to $ \mathds{P}^x$ on $ {\cal X}^0_{t} $, and therefore on $ {\cal X}^0_{t+} $, where $ \mathds{P}^x $ is the law of a standard $\Rd $-valued Brownian motion starting at $x$.
Most of the results of this paper extend straightforwardly to processes that are homeomorphic coordinate changes of the processes described above, for instance to the diffusion process given by the It\^o equation $dX_t = f(X_t)dt + \sigma dW_t, \;\;f$ being an $\Rd$-valued Lipschitz function, $\sigma$ a non-singular $d\times d$-matrix coefficient and $W$ a Brownian motion. Though, for simplicity, such extensions are not generally included in the propositions of this paper, a non-linear coordinate change of this type is employed in the last section for a  generalisation of the work of \cite{fleming1989piecewise} mentioned previously.
 
For a probability measure $\mu$ on ${\cal B}(\Rd)$ let $\mathds{P}^\mu$ denote $\int \mu(dx)\mathds{P}^x$, the law of a Brownian motion $X$ when the distribution of $X_0$ is $\mu$. Generally a superscript $\mu$ (or $x$) on various terms will indicate that at least the distribution of $X_0$ is $\mu$ (or a Dirac measure $\delta_x$ at $x$); thus ${\cal X}^\mu$ is the completion of ${\cal X}^0$ with respect to $\mathds{P}^\mu$ and ${\cal N}^\mu$  the corresponding subset of null sets. ${\cal X}^\mu_{t+}$ denotes the completion of the sub-$\sigma$-algebra ${\cal X}^0_{t+}$ with respect to $\mathds{P}^\mu$ and ${\cal N}^\mu_{t+}$ denotes the corresponding null sets. Note the filtration $\mathbb{X}^\mu = ({\cal X}^\mu_{t+})_{t\geq 0}$, is right continuous\footnote{It is, of course, right continuous in the stronger sense that $\mathbb{X}^\mu = ({\cal X}^\mu_{t})_{t\geq 0}$, but we do not use this fact.} by definition, but it does not satisfy the ``usual conditions" as we are not assuming that ${\cal X}^\mu_{t+}$ contains ${\cal N}^\mu$.
 
The ``initially universal'' sub-$\sigma$-algebra $\bigcap_{\mu}{\cal X}^\mu_{t+}$, where $\mu$ ranges over all probability measures on ${\cal B}(R^d)$,
will be denoted by ${\cal X}^\sim_{t+}$ and the right-continuous filtration $({\cal X}^\sim_{t+})_{t\geq0}$ by $\mathbb{X}^\sim$.
Since its sample paths are continuous, $X$ is progressive,   that is, progressively measurable, with respect to each of the filtrations $\mathbb{X}^o,\;\mathbb{X}^\sim,$ and $\mathbb{X}^\mu$.

Suppose $h$ is a Borel map on $\Rd$ taking values in a separably metrizable\footnote{Here, `separably metrizable' means not only that ${\cal E}$ is countably generated, but also that the single-point sets $\{y\} \in {\cal E}$.} space $\mathbb{E}$ with a Borel field ${\cal E}$.  Let $Y$ denote the $\mathbb{X}^\sim$-progressive $(\mathbb{E}, {\cal E})$-valued process $Y= (Y_t)_{ t \geq 0}$, where $Y_t = h\circ X_t$. Matching the notation for the quantities associated with $X$, ${\cal Y}^\mu_{t+}$ will denote
 $\sigma({\cal Y}^0_{t+}, {\cal N}^\mu_{t+} )$, the augmentation by the null sets of ${\cal X}^\mu_{t+}$ of the natural right-continuous sub-$\sigma$-algebra ${\cal Y}^0_{t+} = \bigcap_{s>t}\sigma(Y_r: 0\leq r \leq s )$ generated by $Y$. Finally ${\cal Y}^\sim_{t+}$ will denote the sub-$\sigma$-algebra $\bigcap_{\mu}{\cal Y}^\mu_{t+}$ and $\mathbb{Y}^\mu$ and $\mathbb{Y}^\sim$ the right-continuous filtrations $({\cal Y}^\mu_{t+})_{t \geq 0}$ and $({\cal Y}^\sim_{t+})_{t \geq 0}$.  $Y$ is intended to represent a running process consisting of a mixture of real measurements and binary observations of $X$. Intuitively, the generated sub-$\sigma$-algebra ${\cal Y}^\sim_{t+}$ represents the information available about $X$, irrespective of its initial distribution, when it is measured through $h$ up to and just after $t$.

The introduction of the filtration $\mathbb{Y}^\sim$, which does not depend on the initial distribution, allows us to replace the global and local problems by a simpler merged version: determine conditions on $h$ under which $X$ can be identified from $Y$ in the sense that $X$ is adapted to $\mathbb{Y}^\sim$.
In several of our results we give conditions under which $X_t$ can be identified for cases where the dimension of the measurement map $h$ is at least as great as that of $X_t$. But this condition on dimension is by no means necessary; the following example shows that there are one-dimensional real=analytic measurement maps that can identify a Brownian motion of any dimension. Its justification is given in the Appendix.

\begin{example}\label{example}
	Suppose $X ;= (X_t)_{t\geq 0}$ is an $\Rd$-valued Brownian motion on the usual canonical space and $h$ is a real-valued function on $\Rd$ of the form
	$ h(x) = e^{a_1x_1} + \ldots + e^{a_dx_d}$, $x\in \Rd$,
	where $a_1, a_2,\ldots, a_d$ are arbitrarily chosen coefficients subject only to the condition that the d-vector $(a_1^{-2},\ldots,a_d^{-2})$ does not lie on any of the $3^d/2$ hyperplanes given by $ \sum_1^d u_jx_j = 0, $ where the $u_j \in \{ -1, 0, 1\}.$ 
	For such an $h$, suppose $Y$ is the continuous real-valued process $Y_t = h(X_t),\; t\geq 0$. Then the Brownian motion $X$ is adapted to the filtration $\mathbb{Y}^\sim$ generated by $Y$.
	
\end{example}

It turns out that for both global and local problems the most natural set of conditions depend on the fine topology of $\mathbb{R}^d$ that arises in classical potential theory. This is the same as the Euclidean topology if $d=1$, but is strictly finer if $d\geq 2$. The features of this topology that are relevant to this paper are outlined below. This summary is based on sections 1.XI.1 p167, 1.XI.6 p177, 2.IX.15 p656 and 2.IX.19 p666 of \cite{doob2012classical}, and sections 3.5 p107 and 4.2 p151 of \cite{chung2006markov}. See also Exercise $(2.25)\; 3^\circ$, p93 of \cite{revuz2013continuous}.

 Suppose $X$ is a Brownian motion with a law from the family of probability measures $\mathbb{P}^x$. For an arbitrary set $A$ in $\mathbb{R}^d$, let $T_A$ be the hitting time
\begin{align*}
T_A = \begin{cases}
\inf\{t>0, X_t\in A\}\\
\infty \quad \text{if } X_t\in A^c \text{ for all } t>0.
\end{cases}
\end{align*}\\
For Borel $A$, $T_A$ is the debut of the progresssively measurable set $\{(t, \omega) \in (0,\infty) \times \Omega : X_t(\omega) \in A\}$ and, as such, is known to be an $\mathbb{X}^x$-stopping time for all $x \in \Rd$.\footnote{A classic theorem that the debuts of progressive sets are stopping times for filtrations that are right-continuous and complete (that is, in the present set-up, each ${\cal X}^0_{t+}$ being augmented with the family of null sets ${\cal N}^x_{t+}$) is given in \cite{revuz2013continuous} as Theorem I.4.15 on p.44. As Revuz and Yor point out, this is a simple consequence of a fundamental theorem from measure theory (see Theorem I.4.14 on the same page) about the measurability of the projections $ (0, t] \times \Omega  \rightarrow \Omega $ for $t>0$ of jointly measurable sets into a probability space.} 
The Blumenthal zero-one law for Brownian motion states that for each $x\in\mathbb{R}^d$ either $\mathbb{P}^x(T_A>0)=1$ or $\mathbb{P}^x(T_A=0)=1$. This dichotomy is exploited in the probabilistic definition of the fine topology: for arbitrary $A$ the set of points $x \in A$ for which
 $\mathbb{P}^x(T_B=0)=1$ for all open $B\supset A$, forms a set $A^r$, the $\mathbf{regular}$ points of $A$. This set corresponds to the set of fine accumulation points of $A$. The $\mathbf{fine \; closure}$ of $A$ is the set $A^* = A\cup A^r$; the $\mathbf{fine \; interior}$ of $A$ is the set  $A^{fo} = A\setminus (A^c)^r, \;\;A^c$ being the complement of $A$. $A$ is finely open if it is its own fine interior and is finely closed if it contains its regular points $A^r$.

For Borel sets $A$, regular and fine interior sets can be defined more simply: $x \in \Rd$ is a regular point of the set $A$ if and only if $\mathbb{P}^x(T_{A}=0)=1$ and $x \in A$ is a fine interior point of $A$ if and only if $\mathbb{P}^x(T_{A^c} > 0) = 1$. Moreover, if $A$ is Borel, its  regular set $A^r$ and its set of fine interior points $A^{fo}$ are also Borel sets (see in Section 4.2, p151 of \cite{chung2006markov}, Fact XII and the subsequent paragraphs). The set $\mathbb{Q}^2 \subset \mathbb{R}^2$, where $\mathbb{Q}$ denotes the rationals, has no regular points, and is an example of a finely closed set that differs markedly from its Euclidean closure $\mathbb{R}^2$.
The non-empty finely open Borel sets $B$ form a basis for the fine topology and have strictly positive $d$-dimensional Lebesgue measures $\lambda_d(B)$. This last statement is a consequence of the contrapositive form of the fact that $\lambda_d$-null sets have complements that are finely dense. Both Euclidean and fine topologies are used in this paper; to distinguish between the two, fine concepts will generally be preceded with the words `fine' or `finely'.
\bigskip

Assume $X$ is a Brownian motion with initial distribution $\mu$. Note that the `observation germ field' ${\cal Y}^\mu_{0+} \subset {\cal X}^\mu_{0+} $, and the inclusion is generally strict. To analyse further the effect of the observation germ field on $X_0$, it is helpful to introduce the following constructs. Let $\widehat{\sigma}^\mu(h)$ denote the class of Borel sets $ \{A \in {\cal B}(\Rd):\: \{X_0 \in A\} \in {\cal Y}_{0+}^\mu \} $, which is clearly a $\sigma$-algebra, and let $\widehat{\sigma}(h)$ denote the intersection $\sigma$-algebra $\bigcap_\mu\widehat{\sigma}^\mu(h)$. We shall call $\widehat{\sigma}(h)$ and $\widehat{\sigma}^\mu(h)$ \textbf{germ enlargements} of $\sigma(h)$. 
The introduction of germ enlargements allows us to `solve' the merged problem mentioned earlier in the sense of rephrasing it purely in terms of $h$:

\begin{lemma}\label{enlargement}
The following are equivalent.
\begin{enumerate}[(i)]
	\item $\widehat{\sigma}(h) = {\cal B}(\Rd)$.\vspace{-5mm}
	\item $X_0$ is ${\cal Y}^\sim_{0+}$-measurable.\vspace{-5mm}
	\item $X$ is adapted to $\mathbb{Y}^\sim$.\vspace{-5mm}
\end{enumerate}

\end{lemma}
\bigskip
\begin{proof}That $(ii)$ implies $(i)$ follows immediately from the definition of the germ enlargements. That $(iii)$ implies $(ii)$ is trivial. So it remains to show that $(i)$ implies $(iii)$. Fix $t>0$. Let $X'$ denote the process $(X_{t+s})_{s\geq 0}$ and $Y'$ the observed process $(h\circ X_{t+s})_{s\geq 0}$. Given the distribution of $X_0$ is $\mu$, denote the distribution of $X_t$ by $\mu_t$. Then the law of $X'$ on its canonical space $\Omega'$ (a copy of $\Omega$) is $\mathds{P}^{\mu_t}. $ So the germ enlargement associated with the observable germ field ${\cal Y'}^{\mu_t}_{0+}$ of $X'$ is $\widehat{\sigma}^{\mu_t}(h)$. 
Since $\widehat{\sigma}^{\mu_t}(h) \supset \widehat{\sigma}(h) = {\cal B}(\Rd),\;X_0'$ is ${\cal Y'}^{\mu_t}_{0+}$-measurable. Consequently $X_t = X'_0 \circ \Theta_t$ is measurable on $\Theta_t^{-1}{\cal Y'}^{\mu_t}_{0+}$, the observable germ field at $t$. Here, $\Theta_t:\; \Omega \mapsto \Omega'$ is the shift operator given by $\Theta_t(\omega)(s) = \omega(t+s),\;s\geq 0$. $(iii)$ then follows from the fact that, for all $\mu$ and $t\geq 0$, ${\cal Y}^\mu_{t+} = \sigma(\Theta_s^{-1}{\cal Y'}^{\mu_s}_{0+},\;0 \leq s \leq t,\; {\cal N}^\mu_{t+}).$ 
 \end{proof}

In line with the notation for generated sub-$\sigma$-algebras of ${\cal B}(\mathbb{R}^d)$, we shall also use $\widehat{\sigma}(h,{\cal C})$ (resp. $\widehat{\sigma}^\mu(h,{\cal C})$) as an alternative notation  for $\widehat{\sigma}(h,\mathds{1}_C :C \in {\cal C})$ (resp. $\widehat{\sigma}^\mu(h,\mathds{1}_C :C \in {\cal C})$), where ${\cal C}$ is some collection of sets in ${\cal B}(\mathbb{R}^d)$.

We have restricted $\widehat{\sigma}(h)$ and $\widehat{\sigma}^\mu(h)$ to sub-$\sigma$-algebras of Borel sets; such sets are sufficiently general for our purposes in that, fortuitously, the fine closure operation maps a Borel set to a Borel set.
 Intuitively, $\widehat{\sigma}(h)$ contains, besides $\sigma(h)$, additional information available about $X_0$ derived from, for instance, the immediate rate of change of any $\mathbb{Y}^\sim$-adapted quadratic-variation process.
 
  To give a specific example where $X_0$ is ${\cal Y}_{0+}^\sim$-measurable, but not $\sigma(Y_0, {\cal N^\sim})$-measurable: suppose $X$ is a Brownian motion in $\mathbb{R}^2$ and $Y_t = h \circ X_t,\; t \geq 0,$ where $h(x_1, x_2) = e^{x_1} - e^{x_2}.$ Then the $\mathbb{Y}^{\sim}$-adapted quadratic variation process of $Y$ is indistinguishable from the continuous process $\int_0^t q(X_s)\mathrm{d}s$, where $q(x_1, x_2) = e^{2x_1} + e^{2x_2}$. Differentiation in $t$ shows that $q(X_0)$ is ${\cal Y}_{0+}^\sim$-measurable; so $q$ is a $\widehat{\sigma}(h)$-measurable function that is clearly not $\sigma(h)$-measurable. However, it is easy to verify in this particular case that the $\mathbb{R}^2$-valued map $(y_1,y_2) \to (\log \frac{(\sqrt{2y_2 - y_1^2}+y_1)}{2},\;\log \frac{(\sqrt{2y_2 - y_1^2}-y_1)}{2})$ on the domain $\{y_2 > y_1^2\}$ is just the inverse of the joint $\widehat{\sigma}(h)$-measurable map $(h,q)$. As a result the $\sigma$-algebras $\widehat{\sigma}(h), \;\sigma(h,q) \; \text{and} \;{\cal B}(\mathbb{R}^2)$ all coincide and $X_0$ is, consequently, ${\cal Y}_{0+}^\sim$-measurable.

The fine topology endows the germ enlargement with some intriguing properties. Some of these, which we shall exploit in later sections, are given in the next two lemmas.

\begin{lemma}\label{germ h}
We have that 
$\widehat{\sigma}(h) = \{ A \in{\cal B}(\Rd):  \mathds{1}_A\circ X \text{ is progressive with respect to } \mathbb{Y}^\sim \}.$ 
Furthermore $\widehat{\sigma}(h) $ is closed under fine closure; that is, if $ A \in \widehat{\sigma}(h)$, then so does its fine closure $A^*$ and its fine interior $A^{fo}$.  
\end{lemma}
 \begin{proof}
Suppose $X$ is a Brownian motion with law $\mathds{P}^\mu$. Then $X$ is $\mathbb{X}^\mu$-progressive by continuity of its paths; that is to say that for each $t \geq 0,$ the $\Rd$-valued  map $(t,\omega) \rightarrow X_t(\omega)$ is measurable on ${\cal B}([0,t])\times {\cal X}^\mu_{t+}$. It follows that $Y = h\circ X$ is not only $\mathbb{X}^\mu$-progressive, but also `self-progressive' in that it is $\mathbb{Y}^\mu$-progressive for all $\mu$ and therefore $\mathbb{Y}^\sim$-progressive.
Suppose that the process $\mathds{1}_A\circ X$ is $\mathbb{Y}^\sim$-progressive for some Borel $A$. This implies that $\{X_0\in A\} \in {\cal Y}^\sim_{0+}$. $A \in \widehat{\sigma}(h)$ then follows from the definition of $\widehat{\sigma}(h)$.
To prove the reverse, choose $ A \in \widehat{\sigma}(h)$. Then $\{ X_0 \in A \} \in {\cal Y}^\mu_{0+}$ for all $\mu$ and $\{ X_t \in A \} \in \Theta_t^{-1}({\cal Y}^{\mu_t}_{0+})$ for $t \geq 0$, where $\Theta_t:\Omega \mapsto \Omega$ is the shift operator $\Theta_t\omega(s)= \omega(t+s),\; s \geq 0$. Since ${\cal Y}_{t+}^{\sim}$ can be expressed as $ \cap_{\mu}\:\sigma(\Theta_s^{-1}({\cal Y}^{\mu_t}_{0+}),0 \leq s \leq t, {\cal N}^\mu_{t+}  )$, $\mathds{1}_A \circ X$ is $\mathbb{Y}^\sim$-progressive. 
 
To prove the second assertion of the lemma, suppose that $ A \in \widehat{\sigma}(h)$. Since $\mathds{1}_A \circ X$ is $\mathbb{Y}^\sim$-progressive, the hitting time  $T_A = \inf\{t >0:\; X_t \in A \}$ is a $\mathbb{Y}^\sim$-stopping time. So the set $\{T_A = 0\}$ is $({\cal Y}_{0+}^\sim)$-measurable. 
Now it follows from Blumenthal's zero-one law that the function $x \rightarrow \mathds{P}^x(T_A = 0)$ takes the values $0 \text{ or } 1$ and from the definition of regular sets that it is just the indicator function $ x \rightarrow \mathds{1}_{A^r}(x)$ of $ A^r$. As previously remarked, $A^r$ can be identified with the set of fine accumulation points of $A$ and is Borel. Since $X$ is Markov we have, for any $\mu$, that

\begin{align*}
\mathds{P}^\mu(X_0 \in A^r,\, T_A > 0) &= \mathds{E}^\mu[\mathds{1}_{\{X_0 \in A^r\}}(1 - \mathds{P}^{X_0 } (T_A = 0))] = \int(\mathds{1}_{A^r}(1 - \mathds{1}_{A^r}) \text{d}\mu = 0,\\
\mathds{P}^\mu(X_0 \in (A^r)^c, T_A = 0) &= \mathds{E}^\mu[\mathds{1}_{\{X_0 \in (A^r)^c\}}\mathds{P}^{X_0 } (T_A = 0)] = \int(1 - \mathds{1}_{A^r})\mathds{1}_{A^r} \text{d}\mu = 0.
\end{align*}

So the symmetric difference $\{X_0 \in A^r\}\bigtriangleup \{ T_A = 0\}$ is a ${\cal Y}^\mu_{0+} $-null set, and $\{X_0 \in A^r\}$ is therefore ${\cal Y}_{0+}^\mu $-measurable, as is also $\{X_0 \in (A^c)^r\}$ by a similar argument based on $T_{A^c}$. Since this holds for all $\mu$, these sets are also ${\cal Y}_{0+}^\sim $-measurable. Furthermore $\{X_0 \in A\}$ is ${\cal Y}_{0+}^\sim $-measurable by definition. As $A^* = A\cup A^r\;\text{and}\; A^{fo} = A \setminus (A^c)^r $, the last assertion of the lemma immediately follows.   
 \end{proof}

A sub-$\sigma$-algebra ${\cal S} \subset {\cal B}(\Rd)$ will be said to be \textbf{separable on a set} $A \in \Rd$,
if the restriction 
${\cal S}|_A = {\cal B}(\Rd)|_A$. Note that this expression implies that not only is ${\cal S}|_A$ countably generated, but one-point sets of $A$ are ${\cal S}|_A$-measurable as well. If $Z$ is an $\Rd$-valued random variable and ${\cal S}$ is generated by a family of measurements $h_\alpha(Z)$ of $Z$, \emph{including $\mathds{1}_A(Z)$}, then the separability of ${\cal S}$ on $A$ would mean that $Z$ is measured exactly whenever its value lies in $A$.

\begin{lemma}\label{dense D}

 Let $D$ be a finely dense Borel set in $\widehat{\sigma}(h)$. Suppose a finely open bounded Borel set $G$ (not necessarily in $\widehat{\sigma}(h)$) has the property that $\widehat{\sigma}(h)$ is  separable on its restriction $G\cap D$. Then $\widehat{\sigma}(h)$ is separable on $G^*$ and its subsets.

\end{lemma}
 
 \begin{proof}
Given that open balls generate ${\cal B}(\Rd)$, the separability of $\widehat{\sigma}(h)$ on ${G\cap D}$ is equivalent to the property that for any open ball $B, \; B\cap G \cap D \in \widehat{\sigma}(h)|_ {G}$. By Lemma \ref{germ h}, $(B\cap G \cap D)^*$ also belongs to $\widehat{\sigma}(h)$.
   But in any topology the closure of the intersection of a dense set with an open set coincides with the closure of the open set. Therefore $ \bar{B} \cap G^* = B^*\cap G^* = (B\cap G)^* = (B\cap G \cap D)^*$, the first equality following from the fact that the Euclidean and fine closures of an open ball $B$ coincide, and the third from the fine density of $D$. 
   Consequently, $\bar{B} \cap {G^*} \in \widehat{\sigma}(h)|_ {G^*}$. Since closed balls generate ${\cal B}(\Rd) $, the restricted $\sigma$-algebras $  {\cal B}(\Rd) |_ {G^*}$ and $ \widehat{\sigma}(h) |_ {G^*}$ coincide, thus establishing the separability of $\widehat{\sigma}(h)$ on $G^*$ and its subsets. This completes the proof.
   \end{proof}
\bigskip

 We need to introduce some topologically localised notions of symmetry and asymmetry. Let  $h:\mathbb{R}^d\to\mathbb{R}^e$ be a Borel map and let $\tau$ be either the Euclidean or fine topology of $\mathbb{R}^d$. Let $G$ be a $\tau$-open Borel set in $\mathbb{R}^d$. For each non-trivial affine isometry $\kappa$ in $\mathbb{R}^d$ that is, an affine map expressible as $\kappa(x)=a + Ox,\; \text{where}\;a$ is a vector and $O$ an orthogonal matrix), let $S_\kappa$ be the $\mathbf{\kappa}$\textbf{-symmetry set} of $h$ in $G$ given by the Borel set $S_\kappa = \{x\in G: \kappa(x)\in G,\;  h\circ\kappa(x) = h(x)\}$. Notice that, as $\kappa|_G$ is a Euclidean affine isomorphism relativized to $G$, it is also a $\tau-$homeomorphism; $S_\kappa$ and its $\tau$-interior map through $\kappa$ onto, respectively, $\kappa S_\kappa = S_{\kappa^{-1}}$ and its $\tau$-interior, both of which are Borel. 
  If $S_\kappa$ has at least a non-empty $\tau$-interior, $h$ will be said to possess a \textbf{partial $\kappa$-symmetry in $G$} with respect to $\tau$.
  The justification for this term is that with any $x \in S_\kappa$ we can associate a bounded neighbourhood $V_x \subset S_\kappa$ that is sufficiently small so that $V_x$ and $\kappa V_x$ are disjoint, and also a piecewise affine map $\widetilde{\kappa}$ that is $\kappa$ on $V_x$, $\kappa^{-1}$ on $\kappa V_x$ and the identity on the complementary set $(V_x \cup \kappa V_x)^c$, such that $h$ is $\widetilde{\kappa}$-invariant in the sense that $h= h\circ\widetilde{\kappa}$.
   If, however, the $\tau$-interior of $S_\kappa$ is empty for all non-trivial affine isometries $\kappa$, $h$ will be said to be \textbf{$\tau$-totally asymmetric in $G$}. If $G$ is $\mathbb{R}^d$ we drop reference to it.    For ${d \geq}2$ fine total asymmetry is generally a stronger (i.e. more specific) property than Euclidean total asymmetry.

\bigskip

We shall also make use of a weaker form of asymmetry that is defined in terms of reflections in hyperplanes; that is, the subclass of affine isometries of the form $\rho(x) = x - 2(\langle n,x\rangle-\alpha)n$ for some unit vector $n$ and real $\alpha$.  
$h$ will be said to possess a \textbf{fine local reflection} if there exists a reflection $\rho$ 
such that the fine interior $S_\rho^{fo}$ of the $\rho$-symmetry set $S_\rho$ contains a fixed point of $\rho$; that is, $S_\rho^{fo}\cap H_\rho $ is nonempty, where $H_\rho$ is the hyperplane of fixed points of $\rho:\; \{x: \langle n,x\rangle = \alpha\}$.
 $h$ will be said to be \textbf{finely totally non-reflective} if there is no such reflection $\rho$. Note that $S_\rho$ contains $H_\rho$ and, since $\rho^2$ = id, both $S_\rho$ and $S_\rho^{fo}$ are $\rho$-invariant.
We shall need a property of fine local reflections which is a consequence of a slight variation of a theorem of Gardner (theorem A in [G]) which itself is an enhancement of a theorem of Lyons (theorem 4.4 in \cite{lyons1980finely}) on the polygonal connectedness of finely connected finely open sets. Let $[x,y]$ denote the closed line segment between $x$ and $y$ in $\mathbb{R}^d$, and $p_{[x,y]}$ the hyperplane $\{z\in\mathbb{R}^d: |z-x| = |z-y|\}$ bisecting it. 

\begin{lemma}[Gardner \& Lyons]\label{GardnerLyons}
Let $V$ be a finely open set in $\mathbb{R}^d$ and let $x_0\in V$. Then there is a fine neighbourhood $V_0\subset V$ of $x_0$ with the property that, for any distinct points $x,y\in V_0$, there is a set $A\subset p_{[x,y]}$ of strictly positive $(d-1)$-dimensional Lebesgue measure, such that 
\begin{equation}\label{gl}
|z-y|=|z-x| < \frac{|x-y|}{\sqrt{3}} \text{ and } [x,z]\cup [z,y]\in V \text{ for all } z\in A.
\end{equation}
\end{lemma}

The above statement also holds if we replace the constant $\sqrt{3}$ in (\ref{gl}) with any other positive constant $a>1$.     

A simple but useful consequence of this lemma appropriate for the case where $V$ is invariant with respect to a reflection is as follows:

\begin{lemma}\label{CorrolG-L} 
 If $\rho$ generates a fine local reflection for a Borel map $h$, $H_\rho$ is its hyperplane of fixed points and $S_\rho$ is its symmetry set for $h$, then there is a finely open simply connected $\rho$-invariant Borel set $V_1\subset S_\rho$ such that its intersection $V_1\cap H_\rho$ with $H_\rho$ has a strictly positive $(d-1)$-dimensional Lebesgue measure.
\end{lemma}
\begin{proof} First note that $S_\rho$ is Borel. Let $x_0 \in S_\rho^{fo}\cap H_\rho $; that is, $x_0$ is a fixed point of $\rho$ in the $\rho$-invariant fine interior  $S_\rho^{fo}$ of its symmetry set. Take $V$ of Lemma \ref{GardnerLyons} to be the Borel set $S_\rho^{fo}$. If $V_0$ is a Borel neighbourhood of $x_0$ with the properties of that lemma then $V_1 = V_0 \cap \rho V_0$ is a $\rho$-invariant neighbourhood with the same properties. Pick $x_1$ in the $\rho$-invariant set $V_1\setminus H_\rho$. Then $\rho x_1$ is distinct from $x_1$ and $H_\rho$  coincides with the hyperplane $p_{[x_1,\rho x_1]}$ of Lemma \ref{GardnerLyons}; the present lemma immediately follows.
\end{proof}

\section{Necessary Conditions}
The next two propositions illustrate how some form of asymmetry for $h$ that involves the fine topology is required for the solution of the global and local problems.

\begin{proposition}\label{globalnecessity}  

Let $\mu$ denote the distribution of $X_0$. We make the following assumptions:

\begin{enumerate}[(i)]
	\item $h$ is a Borel function\vspace{-5mm}
	\item the law $\mathds{P}$ of $X$ is locally equivalent to the law ${\mathds{P}}^{\mu}$ of a Brownian motion (in the sense that for \vspace{-3mm}
	
	each finite $t\geq 0$ both probability measures induce the same completion ${\cal X}^\mu_{t+}$ of ${\cal X}^0_{t+}$)\vspace{-5mm}
	\item $X_0$ is ${\cal Y}^\mu_{0+}$-measurable\vspace{-5mm}
	\item $G$ is a non-empty finely open Borel set on which $\mu$ has a strictly positive probability density (that \vspace{-10mm}
	
	is, the restricted measure $\mu|_G$ is equivalent to $\lambda_d|_G$).\vspace{-5mm}
		
\end{enumerate}	

\noindent Then $h$ is finely totally asymmetric on $G$. In particular, if $\mu$ has a strictly positive probability density everywhere and $X_0$ is ${\cal Y}^\mu_{0+}$-measurable, then $h$ is  finely totally asymmetric.

\end{proposition}

The proof relies on establishing, contrapositively, that, if $h$ possesses a local $\kappa$-symmetry in $G$, then there exists a bounded Borel finely open subset $B$ of $G$ such that its disjoint $\kappa$ image $\kappa B \subset G$ and $h|_B = h\circ \kappa|_B$. Intuitively, A path that starts in $B$ will remain in it for some time and will be `confused' with a path starting in $\kappa B$.

\begin{proof}
Clearly we may take $X$ to be a Brownian motion with law $\mathds{P}^\mu$. Suppose, for some affine isometry $\kappa, \;h$ possesses a local $\kappa$-symmetry for some finely open Borel $G$ with the properties of the proposition; that is, the Borel symmetry sets $S_\kappa $ and $\kappa (S_\kappa)$ in $G$ have non-empty finely open interiors, the non-empty finely open Borel subsets of which are all of positive $\lambda_d$ measure and therefore of positive $\mu$ measure.
First note that the affine sub-space $H_\kappa$ of fixed points of $\kappa$, of dimension at most $d-1$, is a closed set, and therefore a finely closed set. Furthermore, $\lambda_d(H_\kappa) = 0$; so the finely open sets $S^{fo}_\kappa \setminus H_\kappa$ and $\kappa S_\kappa^{fo} \setminus H_\kappa$ are non-empty and possess positive $\mu$ measures.

 Now choose a point $x$ in $S^{fo}_\kappa \setminus H_\kappa$ . Since the fine topology is Hausdorff and open balls centred at $x$ are also fine neighbourhoods of $x$ we can choose a finely open neighbourhood $B$ of $x$ within $S^{fo}_\kappa \setminus H_\kappa$ that is Borel and bounded, and disjoint from its bounded Borel image $\kappa B$ in $\kappa S^{fo}_\kappa\setminus H_\kappa$.

Let $\widetilde{\kappa}:\mathbb{R}^d \to \mathbb{R}^d$ denote the reflexive Borel isomorphism that is $\kappa$ on $B$, $\kappa^{-1}$ on $\kappa B$ and the identity elsewhere. 
Since $B\cup \kappa B\subset S_\kappa \cup \kappa S_\kappa $, $h$ is locally $\widetilde{\kappa}$-invariant in the sense that the restricted function $h|_{B\cup \kappa B}$ is measurable with respect to the $\sigma$-algebra ${\cal I}(\widetilde{\kappa}) = \{ A \in {\cal B}(\Rd): A= \widetilde{\kappa}A\}$ of Borel sets in $\mathbb{R}^d$ that are invariant under $\widetilde{\kappa}$. It is useful to note that ${\cal I}(\widetilde{\kappa})$ is also $\sigma(\kappa_B)$, the $\sigma$-algebra generated by the $\Rd $-valued map $ \kappa_B = \kappa \mathds{1}_B + \text{id}\:\mathds{1}_{B^c}$ where id is the identity on $\Rd$. As such it is countably generated and has set atoms that are either the pairs of the form $\{  x, \kappa x \}$ that partition $ B\cup \kappa B $ or the singletons $\{x\}$ that make up $(B\cup \kappa B)^c $.

Since processes with locally equivalent laws share the same locally completed filtrations we may take the law of $X$ to be $\mathbb{P}^{\nu}$ where $\nu$ is any probability measure equivalent to $\mu$. Note that $\mu$ restricted to the set $B \cup \kappa B \subset G$ is equivalent to Lebesgue measure on that set. A suitable choice of $\nu$ can be  defined as follows.  Denote  the bounded set $B \cup \kappa B$ by $ \widetilde{B}$ for short. For Borel $D$ in $\mathbb{R}^d$ set
\begin{align*}
\nu(D) = \frac{\mu(\widetilde{B})}{\lambda_d(\widetilde{B})}(\lambda_d(D\cap \widetilde{B}) + \mu(D\setminus \widetilde{B}) .
 \end{align*}
Given that Lebesgue measure is invariant for the group of affine isometries acting on $\Rd$, it is clear that $\nu$ is $\widetilde{\kappa}$-invariant in the sense that $\nu(D) = \nu(\widetilde{\kappa}^{-1}D)$.

The map $\widetilde{\kappa}$ on $(\mathbb{R}^d,{\cal B}(\mathbb{R}^d) ,\nu)$ induces a map with similar invariance properties on $(\Omega, {\cal X}^\nu, \mathds{P}^\nu)$ in the following way. First note that $(\Omega, {\cal X}^0)$ is separable (it is metrizable, for instance, with any metric giving point-wise convergence of sequences $(\omega_n)$ restricted to rational $t$).
Let $\widetilde{K}$ be the map that for $t \geq 0$ is given by

\begin{align*}
\widetilde{K}_t(\omega) = \begin{cases}
\kappa\circ X_t(\omega), \qquad & \text{ if } X_0(\omega)\in B,\\
\kappa^{-1}\circ X_t(\omega), \qquad & \text{ if } X_0(\omega) \in \kappa B,\\
X_t(\omega), \qquad & \text{ if } X_0(\omega) \in (B\cup \kappa B)^c.
\end{cases}
\end{align*}
Then $\widetilde{K}$ determines a  Borel isomorphism on $(\Omega, {\cal X}^0)$ onto itself.

The spatial homogeneity of Brownian motion means that the laws $\mathds{P}^x$ possess the property that $\mathds{P}^x(\widetilde{K}^{-1}(A)) = \mathds{P}^{\widetilde{\kappa}(x)}(A)$ for ${\cal X}^0$-measurable $A$. This property, the $\widetilde{\kappa}$-invariant property of $\nu$ and the disintegration of $\mathds{P}^\nu$ as  $\int \nu(dx)\mathds{P}^x$, taken together, imply that $\mathds{P}^\nu$ is $\widetilde{K}$-invariant. 
Let $T$ be the first hitting time $\inf \{ t> 0, X_t\in (B\cup \kappa B)^c\}$. Then for each $x \in B\cup \kappa B, \;\mathds{P}^x(T>0) =1$. Consequently, since $\nu$ has a positive density on the finely open set $B\cup \kappa B,\;\nu(B\cup \kappa B) >0$. Now consider the set $\{T>t\}$. This splits into disjoint sets: 
\begin{align*}
\{T>t\} =  B_{(0,t]}\cup B_{(0,t]}',
\end{align*}
where $B_{(0,t]} = \{X_s \in B, s\in (0,t] \}$ and $B_{(0,t]}' = \{X_s \in \kappa B, s\in (0,t] \}$.  

It follows from the definition of $\widetilde{K}$ that $\{T>t\}$ is $\widetilde{K}$-invariant.

Now let $J$ be the $\sigma(X_0)$-measurable variable
\begin{equation*}
J = \mathds{1}_{\{X_0\in B\}} - \mathds{1}_{\{X_0\in \kappa B\}} 
\end{equation*}
Then $\mathds{E}^\nu[J] = 0$. Since $\mathds{E}^\nu[|J|] = \mathds{P}^\nu(X_0\in B\cup \kappa B)>0$, $J$ is nontrivially random. So to prove that ${\cal Y}^\nu_{0+}$ is strictly smaller than ${\cal X}^\nu_{0+}$ whenever $h$ has a local symmetry it is sufficient to show that $J$ differs from its projection $\mathds{E}^\nu[J|{\cal Y}^\nu_{0+}]$ with positive $\mathds{P}^\nu$ probability.  Consider the limit in probability 
\begin{equation*}
\mathds{E}^\nu[J|{\cal Y}^\nu_{0+}] =\text{p}\lim_{t\downarrow 0}\mathds{E}^\nu[J|{\cal Y}^\nu_{t}].
\end{equation*}
It is enough  to show that this limit exists and vanishes. Expand $\mathds{E}^\nu[J|{\cal Y}^\nu_{t}]$ as follows:
\begin{equation*}
\mathds{E}^\nu[J|{\cal Y}^\nu_{t}] = \mathds{E}^\nu[J\mathds{1}_{\{T=0\}}|{\cal Y}^\nu_{t}] + \mathds{E}^\nu[J\mathds{1}_{\{0<T \leq t\}}|{\cal Y}^\nu_{t}] + \mathds{E}^\nu[J\mathds{1}_{\{T>t\}}|{\cal Y}^\nu_{t}]
\end{equation*}
Since $(B\cup \kappa B)^c$ is finely closed, $\{T=0\} = \{X_0\subset (B\cup \kappa B)^c\}$, on which $J=0$, implying that the first term in the expansion is zero. The second term satisfies 
\begin{align*}
\mathds{E}^\nu\left[\Big|\mathds{E}^\nu[J\mathds{1}_{\{0<T \leq t\}}|{\cal Y}^\nu_{t}]\Big|\right] \leq \mathds{P}^\nu(0<T\leq t)
\end{align*}
and vanishes in first moment as $t\to 0$. So the proof is complete if we can show that the third term in the expansion is zero for all $t>0$.  Let ${\cal C}_t$ be the $\pi$-system of ${\cal Y}^\nu_{t}$-measurable cylinder sets of the form $\{h\circ X_{t_i} \in A_i, i = 1,\hdots, n\}$ where $n\in\mathbb{N}$, $t_1,\hdots,t_n\in[0,t]$ and $A_1,\hdots,A_n\in B(\mathbb{R}^d)$. Along with ${\cal N}^\nu_t$, the null sets of ${\cal X}^\nu_t$, ${\cal C}_t$ generates ${\cal Y}^\nu_{t}$; so we only need to show $\mathds{E}[J\mathds{1}_{\{T>t\}\cap C}]=0$ for all $C\in{\cal C}_t$. 
However, any $C\in{\cal C}_t$ is locally $\widetilde{K}$-invariant in the sense that its intersection with $\{T>t\}$ can be written as 
\begin{align*}
C\cap(B_{(0,t]}\cup B_{(0,t]}') = \{ X_{t_i}\in I_i, i = 1,\hdots,n\}\cap(B_{(0,t]}\cup B_{(0,t]}')
\end{align*}
where $I_1,\hdots, I_n\in{\cal I}(\widetilde{\kappa})$. But $\mathds{E}^\nu[J\mathds{1}_{\{T>t\}\cap C}] = \mathds{P}^\nu(C\cap B_{(0,t]}) - \mathds{P}^\nu(C\cap B_{(0,t]}') = 0$ by the $\widetilde{K}$-invariance of $\mathds{P}^\nu$ and the fact that $C\cap B_{(0,t]}$ and $C\cap B_{(0,t]}'$ are images of each other under $\widetilde{K}$. This completes the proof of the first assertion. Since $\Rd$ is trivially the fine interior of the fine support of any $\mu$ equivalent to Lebesgue measure, the second assertion is just a special case of the first where $G$ is $\Rd$.
\end{proof}
\bigskip

The corresponding proposition for the local problem of continuation requires some preliminary explanation. Suppose the measurement function $h$ possesses a a local $\kappa$-symmetry for some $\kappa$ with at least one fixed point that generates a finite cyclic group of $k$ isometries. Suppose a Brownian motion $X$ with initial distribution $\mu$ encounters a fixed point of such a $\kappa$  within the fine interior of its symmetry set. Then, intuitively, the local problem appears to be ill-posed in that the $\mathbb{Y}^\mu$-adapted process of the conditional law of $X_t$, even if it has a left-continuous limit that is a unit point-mass (Dirac) measure $\delta_{ X_{t-} }$, could split, after the encounter, into a member of a family of multiple discrete measure processes supported by a $\mathbb{Y}^\mu$-adapted set process of $k$ points. So, when considering the local problem of continuation and the properties of $h$ that are required for $ X$ to be adapted to $\mathbb{Y}^x$, one expects such affine isometries with fixed points to play a role. Distinct among these, however, are the reflections, with their $(d-1)$-dimensional affine subspaces of fixed points. The affine subspaces of the fixed points of all other cyclic affine isometries have dimensions $d-2$ or less and have the property that a Brownian motion that does not start at a fixed point of the isometry will, with probability one, never meet one. Consequently the presence of a local symmetry set of $h$ for such an isometry has an insignificant effect. In contrast, a Brownian motion will eventually meet any finely open part of the hyperplane of fixed points of a local reflection with positive probability; after the random time when this occurs, $X$ will no longer be adapted to $\mathbb{Y}^x$.

The condition of totally fine asymmetry is no longer required for a solution of the local problem. However, it turns out that the weaker condition of total non-reflectivity of $h$ does remain a requirement, as the following proposition shows.

\begin{proposition}\label{localnecessity}
 Suppose a Brownian motion  $X$ with initial distribution $\mu$ is adapted to the filtration $ \mathbb{Y}^{\mu}$ generated by $h$. Then $h$ is totally non-reflective with respect to the fine topology.
\end{proposition}
 
 The following lemma is useful in the proof of this proposition:

 \begin{lemma}\label{andre}
 Suppose $h$ possesses a fine local reflection $\rho$ with a symmetry set $S_\rho$ and the hyperplane of fixed points of the  reflection is described by $H_\rho =\{ x: x^1 = 0\}$. Suppose the initial distribution $\mu$ is supported by $H_\rho$ and satisfies  $ \mu(S_\rho^{fo}\cap H_\rho) > 0 $. Then $\mathbb{X}^\mu \neq \mathbb{Y}^\mu$.
 \end{lemma}
 
 \begin{proof}[Proof of Lemma \ref{andre}.]
Note that $X^1_0 =0$. By the definition of local reflections we can choose a finely open $\rho$-invariant bounded Borel set $V$ that is contained in the fine interior $ S_\rho^{fo}$ of $S_\rho$ such that $\mu(V\cap H_\rho)>0$. Let $T$ be the hitting time of the complementary set $V^c$. Then $\mathds{P}^z(T > 0) > 0$. Since for $t\geq 0, \; \{T > t\} = \{(|X^1_s|,X^2_s, \cdots,X^d_s) \in V :0 < s \leq t\}$, $T$ is a $\mathbb{Y}^\mu$-stopping time. We now construct an alternative observation process $Z^T $ that is simpler and `more informative' than $Y$. We do this in two steps. First, let $g: \Rd \rightarrow \Rd$ denote the map

 $g(x) = (|x^1|, x^2,\cdots, x^d)\mathds{1}_{V}(x) + \text{id}\:\mathds{1}_{V^c}(x)$
 
\noindent and introduce the $\Rd$-valued process $Z_t = g \circ X_t,\;t\geq 0$. Let ${\cal Z}^\mu_t$ denote the generated sub-$\sigma$-algebra  $\sigma(Z_s: 0\leq s \leq t, {\cal N}^\mu_t)$. Clearly $\sigma(h) \subset \sigma(g)$; so the filtration $({\cal Z}^\mu_{t+})$ is `more informative' than $\mathbb{Y}^\mu$ in the sense that ${\cal Z}^\mu_{t+} \supset {\cal Y}^\mu_{t+}$ for $t \geq 0$. $A\;fortiori,\; T$ is a $({\cal Z}^\mu_{t+})$-stopping time. Secondly, let $Z^T$ be the enhanced observation process 

$Z^T_t = Z_t \mathds{1}_{\{t<T\}} + X_t \mathds{1}_{\{t\geq T\}}, \quad t \geq 0$.

\noindent The corresponding right-continuous generated filtration $({\cal Z}^{T\mu}_{t+})$ turns out to have the simple `interpolated' form \footnote{This can be thought of as a  filtration interpolated between  the filtration $ ({\cal Z}^{\mu}_{t+}) $ and the greater filtration $ {(\cal X}^{\mu}_{t+})$. This is closely related to a construction suggested by M.Yor for a progressively expanded filtration that incorporates a random variable as a stopping time.
 See \cite{yor1978grossissement} and \cite{protter2003stochastic} p.370. }

${\cal Z}^{T\mu}_{t+}  = \{A\in {\cal X}^\mu_{t+} : A \cap \{t< T \} \in {\cal Z}^\mu_{t+} \}, \quad t \geq 0$.

\noindent Furthermore, for $t\geq 0,$

 ${\cal Y}^\mu_{t} \subset {\cal Z}^\mu_{t} \subset {\cal Z}^{T\mu}_{t} \subset {\cal X}^\mu_{t} \text{ and } {\cal Y}^\mu_{t+} \subset {\cal Z}^\mu_{t+} \subset {\cal Z}^{T\mu}_{t+} \subset {\cal X}^\mu_{t+}$.

More precisely, let ${\cal C}_t$  be the $\pi$-system of cylinder subsets $A \in \Omega$ of the form $\bigcap^n_{j=1}\{Z^T_{t_j} \in B_j \}$ where the $t_j \in [0,\infty) \text{ with } t_1 < t_2< \cdots < t_n \leq t$ and the $B_j \in {\cal B}(\mathbb{R}^{d})$. 
Then ${\cal C}_t$ generates (after augmentation with ${\cal N}^{\mu}_t $-null sets) ${\cal Z}^{T\mu}_t $ and the collection $({\cal C}_t)$ generates (after the process of right continuity) the right-continuous filtration $({\cal Z}^{T\mu}_{t+}). $
Note, however, that each $A \in {\cal C}_t $ can also be expressed as the disjoint union $(\bigcap^n_{j=1}\{X_{t_j} \in B_j\cup \rho B_j ,\:t < T \} )\cup (\bigcap^n_{j=1}\{X_{t_j} \in B_j ,\:t \geq T \})$.
Choose $r,t,\;0 <r <t$ so that $\mathds{P}^\mu (|X_r|>0,\:t<T) >0$.
 The $\rho$-invariance of the joint distributions of the Brownian vector variables $X_{t_1}, X_{t_2}, \cdots, X_{t_n}$ when restricted to $\{t < T\}$ show that for each $A \in {\cal C}_t $ 

$\mathds{P}^\mu(\{ X^1_r >0, t<T\}\cap A) = \mathds{P}^\mu(\{ X^1_r <0,t<T\}\cap A)$
 
\noindent and that the non-trivial ${\cal X}^\mu_{t+}$-measurable variable

$J =  \mathds{1}_{\{X^1_r >0,\:t<T\}} - \mathds{1}_{\{X^1_r < 0,\:t<T\}}$

\noindent has zero mean when conditioned on ${\cal Z}^{T\mu}_{t+}$ and therefore when conditioned on ${\cal Y}^\mu_{t+}$. Thus ${\cal Y}^\mu_{t+}$ is strictly smaller than ${\cal X}^\mu_{t+}$. This completes the proof of the lemma.
\end{proof} 
\begin{proof}[Proof of Proposition \ref{localnecessity}.] It suffices to prove the contrapositive assertion that the existence of a fine local reflection $\rho$ for $h$ implies that the filtration $\mathbb{Y}^\mu$ is strictly smaller than $\mathbb{X}^\mu$. If, with the notation of Lemma \ref{andre}, $\mu (S_\rho^{fo}\cap H_\rho) > 0$, then the proposition would follow immediately from the lemma.
  
   For the general case, let $R$ is the first entry time not less than 1 of the hyperplane $H_\rho$ of fixed points of $\rho$. Then $\mathds{P}^\mu (X_R \in H_\rho) = 1$. We assert that the induced distribution of $X_R$ on $H_\rho$ is equivalent to $\lambda_{d-1}, \;(d-1)$-dimensional Lebesgue measure on $H_\rho$. To justify this, note that $X_1$ has a density, that $\mathds{P}^\mu (X_1 \in H_\rho) = 0$, that the induced distribution 
   
   $\mu_R(H_\rho\cap A) := \mathds{P}^\mu (X_R \in H_\rho \cap A) =\mathds{E}^\mu[ \mathds{P}^{X_1} (X_{R-1} \in H_\rho\cap A)], \text{ for } A\in {\cal B}(\Rd)$
   
\noindent by the Markov property and that $\mathds{P}^{x} (X_{R_0} \in H_\rho\cap A)$, where $R_0$ is the first entry time of $H_\rho$, has a Cauchy distribution \cite{revuz2013continuous} with a strictly positive density for all $x \in H^c_\rho$. Consequently by Lemma \ref{CorrolG-L} $ \mu_R(S_\rho^{fo}\cap H_\rho) > 0 $. Now introduce the observation process $ Y':\:Y'_t = X_t \:\mathds{1}_{\{R \geq t \geq  0\}} + Y_t \: \mathds{1}_{\{R < t \}}$, which is clearly more informative than $Y$. Since $X$ is strong Markov and $Y'_R = X_R$ it is enough to show that the delayed Brownian motion $(X_{R + s})_{s\geq 0}$, with initial distribution $\mu_R$, is not observed exactly by the delayed process $(Y'_{R + s})_{s\geq 0}$. This follows from Lemma \ref{localnecessity}; so the proof is complete. \end{proof}

\section{Sufficient Conditions}
This section is concerned with conditions on $h$ and the initial distribution $\mu$ that ensure that $X_0$ is at least ${\cal Y}^\sim_{0+}$-measurable and $X$ is adapted to the filtration $ \mathbb{Y}^\sim$.
 We first consider a subclass of measurement functions for which the condition of fine total asymmetry suffices. In the following ${\cal Y}^\sim_{t,t+}$ will denote the``germ field" $\Theta^{-1}_t{\cal Y}^\sim_{0+}$, which is a sub-$\sigma$-algebra of ${\cal Y}^\sim_{t+}$ generated by $Y$ at and `just after' $t$. See \cite{knight1970remark} for more on germ fields.

\begin{proposition}\label{incongruence}
Suppose that the Borel map $h:\mathbb{R}^d \to \mathbb{R}^e, e \geq d,$ is continuously differentiable on a Euclidean open set $G$ that is finely dense in $\Rd$, is $\widehat{\sigma}(h)$-measurable and supports a  Jacobian matrix $Dh:=\left[\frac{\partial h^i}{\partial x_j}\right]$ of full rank $d$. If $d \geq 2$, further suppose that $h$ is twice continuously differentiable on such a $G$.
 If $h$ is totally asymmetric with respect to the fine topology, then
\begin{enumerate}[(i)]
\item $X_0$ is measurable on ${\cal Y}^\sim_{0+}$.
\item $X$ is adapted to the filtration $\mathbb{Y}^\sim$.
\end{enumerate}
\end{proposition}

The following is a simplification that serves for most purposes.  It also provides examples in which the crucial open set $G$ is both finely dense and measurable with respect to $\widehat{\sigma}(h)$ rather than to $\sigma(h)$.
\begin{corollary}\label{simplification}
Suppose $h:\Rd \to \mathbb{R}^e, \;e\geq d$ is twice continuously differentiable and the Jacobian matrix $Dh$ is of full rank $d$ except on a Lebesgue-null set. Then $X$ is adapted to $\mathbb{Y}^\sim$ if and only if $h$ is finely totally asymmetric.
\end{corollary}

{
\begin{proof}[Proof of the corollary]
  \textit{The `if' part.} Since $h$ is twice continuously differentiable  the process $Y=h \circ X$ is for every $\mu$ a $\mathbb{Y}^\mu$-adapted semimartingale that engenders a matrix covariation process with continuously differentiable paths which, by classic theorems on the approximation of the quadratic variation of a semimartingale (see e.g. \cite{protter2003stochastic} Ch.II), is similarly adapted. Consequently the random matrix of derivatives $Dh(X_0)Dh(X_0)^T$ is ${\cal Y}^\mu_{0+}$-measurable for every $\mu$ and $DhDh^T$ is $\widehat{\sigma}(h)$-measurable. Take $G$ to be the open $\widehat{\sigma}(h)$-measurable set $\{x: \text{det}(Dh^TDh(x))>0\}$ on which $Dh$ is of full rank. Since we are assuming $\lambda_d(G^c) = 0$, $G^c$ has no fine interior. So $G$ is a finely dense set in $\Rd$ and the `if' part of the corollary then follows from the proposition. \textit{The `only if' part.} If $X$ is adapted to $\mathbb{Y}^\sim$ then it is adapted to $\mathbb{Y}^\mu$ for any $\mu$ with a strictly positive density. It then follows from Proposition \ref{globalnecessity} that $h$ is finely totally asymmetric.
\end{proof}

\begin{proof}[Proof of Proposition \ref{incongruence}]

The first part of the proof is mainly analytical in nature. The strategy is to exploit the inverse function theorem to generate a countable number of $(\mathbb{R}\cup \{\infty\})^d$-valued $\sigma(h,G)$-measurable functions $u_n$, each of which is the identity on an open ball, so that for any $t\geq 0$ at least one of the ${\cal Y}^\sim_{t,t+}$-measurable variables $u_n(X_t)$ coincides with $X_t$. Here, $\sigma(h,G)$ is the Borel sub-$\sigma$-algebra generated by $h$ and $G$.  In the second, more probabilistic, part of the proof additional $({\cal Y}^\sim_{t,t+})$-adapted `test' processes are constructed that determine which of the `candidate' processes $u_n(X_t)$ exactly match $X_t$ at any particular time $t$.

The functions $u_n$ are defined on the elements of a countable covering of $G$ by open balls. Let $B(x,r)$ denote an open ball with centre $c$ and radius $r$ and $\bar{B}(x,r)$ its closure. Consider the family ${\cal C}$ of balls of the form $\{ B(x,r)\subset G:x\in G  , 0< r \leq \frac{1}{2} (d(x,G^c)\land 1)\}, \;d(x,G^c)$ being the distance between $x$ and $G^c$. Then ${\cal C}$ covers $G$. Since $\Rd$ is $\sigma$-compact and locally compact and the rank of $ Dh = d$ on $G$, the inverse function theorem allows us to select sequences of points $\{c_1,c_2,\hdots\}$ in $G$ and radii $\{r_1, r_2, \hdots\}$ so that $\{B(c_n,r_n)\}_n $ is an open sub-covering in ${\cal C}$ with the following properties:

\begin{enumerate}[(a)]
\item $\bigcup_n B(c_n,r_n) = \bigcup_n \bar{B}(c_n,r_n) = G$, 


\item for each $n$, $h$ is a diffeomorphism on $\bar{B}(c_n,r_n)$ of class ${\cal C}^1$ if $d=1$ and of class ${\cal C}^2$ if $d\geq 2$.
\end{enumerate}

Let $B_n$ denote $B(c_n,r_n)$. By assumption $h$ is continuous on $G$ and by (b) $h(B_n)$ is open. Let $A_n$ denote the open set $ \{x\in G: h(x)\in h(B_n\}$.
 Let $h^{-1}_n: h(B_n)\to B_n$ denote the  ${\cal C}^1$ or ${\cal C}^2$ diffeomorphism that is the inverse of the restriction $ h_n = h|_{B_n}$. 
 Note that, for any $m$ and $n$, if $x \in B_m \cap B_n$ then $h_m^{-1}\circ h(x) = h_n^{-1}\circ h(x)$.
 
 Set $u_n$ to be the $B_n \cup \{0^d \}$-valued map on $\Rd$ ($0^d$ being the origin in $\Rd$) given by

\begin{align*}
u_n = \begin{cases}
h^{-1}_n\circ h &\quad \text{on }\, A_n\\
0^d &\quad \text{on }\, A^c_n
\end{cases}
\end{align*}

\noindent Then the restricted map $u_n|_{A_n}$ is of class ${\cal C}^1$ or ${\cal C}^2$ and its image is $B_n$; furthermore $u_n$ is a local identity in the sense that $u_n|_{B_n} = \text{id}|_{B_n} $, where id is the identity map. In addition, since $A_n = G\cap h^{-1}( h(B_n))$,   $\; A_n \in \sigma (h,G)$ and $u_n$ is $\sigma (h,G)$-measurable. The countable family of $\Rd$-valued functions $u_n$ can be regarded as ``substitute'' measurements that generate the $\sigma$-algebra $\sigma(h|_G,G)$. But they also serve to expose the presence of `aliases', where $h(B_n)\cap h(B_m) $ is non-empty for some disjoint $B_n$ and $B_m$. Note that, \emph{a priori, $B_n$ is not necessarily in $\sigma (h,G)$}. A crucial part of the proof will be to show that, with the assumption of fine total asymmetry, $B_n \in \widehat{\sigma} (h,G)$.
 
 As will be explained in the second part of the proof, the following functions are needed in the construction of `test' processes that we use to determine which of the processes $u_n(X_t)$ coincide with $X_t$ at any time. Consider the family
\begin{align*}
v_n^{ij} = \begin{cases}
\sum_l \dfrac{\partial u_n^i}{\partial x^l}\dfrac{\partial u_n^j}{\partial x^l} &\quad \text{on }\, A_n\\
0 &\quad \text{on }\, A^c_n
\end{cases}
\end{align*}
for $i,j\in \{1,\hdots,d\}$ and $n\in\mathbb{N}$. Such functions are continuously differentiable on $A_n$ (just continuous if $d=1$). For $d\geq 2$, let
\begin{align*}
w_n^{ij,k} = \begin{cases}
\sum_m \dfrac{\partial v_n^{ij}}{\partial x^m}\dfrac{\partial u_n^k}{\partial x^m} &\quad \text{on }\, A_n\\
0 &\quad \text{on }\, A^c_n
\end{cases}
\end{align*}
for $i,j,k\in \{1,\hdots,d\}$ and $n\in\mathbb{N}$. It will turn out that all the processes $v_n^{ij}(X_t)$ and $w_n^{ij,k}(X_t)$ are ${\cal Y}^\sim_{t,t+}$-measurable. Notice that $w_n^{ij,k}$ can be expanded as 
\begin{align}\label{chris}
w_n^{ij,k} = q_n^{ik,j} + q_n^{jk,i}
\end{align}
for $i,j,k\in \{1,\hdots,d\}$ and $n\in\mathbb{N}$, where
\begin{align*}
q_n^{ij,k} = \begin{cases}
\sum_{l,m}\frac{\partial^2 u_n^k}{\partial x^l \partial x^m}\frac{\partial u_n^i}{\partial x^l}\frac{\partial u_n^j}{\partial x^m} &\quad \text{on }\, A_n\\
0 &\quad \text{on }\, A^c_n.
\end{cases}
\end{align*}
A simple rearrangement of the family of equations (\ref{chris}) shows that 
\begin{align*}
q_n^{ij,k} = \frac{1}{2}(w_n^{ik,j} - w_n^{ij,k} + w_n^{jk,i}).
\end{align*}
We note without further comment that this mimics a familiar rearrangement, dating back to Riemann, for the inversion of the triply ordered family of equations relating the gradients of the components of a Riemannian metric to its Christoffel symbols of the first kind. (See, for instance, Spivak \cite{spivak1979comprehensive} Vol.2 pp.185-186).\\

The following two lemmas about the functions $u_n$ and their first and second order partial derivatives will be needed in the proof of Proposition \ref{incongruence}. They characterize the property of total asymmetry in terms of these derivatives.
\begin{lemma}\label{coincidence}
Suppose $d\geq 2$. For each $n$, the set in $A_n$, on which the Jacobian matrix $Du_n$ is orthogonal and the second-order partial derivatives of $u_n$ are all zero, given by
\begin{align*}
\widetilde{Z}_n = \{x \in A_n: \quad \sum_l \dfrac{\partial u_n^i}{\partial x^l}\dfrac{\partial u_n^j}{\partial x^l}(x)=\delta_{ij},\quad \frac{\partial^2 u_n^k}{\partial x^i \partial x^j}(x) = 0,\quad i,j,k\in \{1,\hdots,d\}\},
\end{align*}
contains $B_n$ and coincides with the set
\begin{align*}
Z_n = \{x \in A_n:\quad v_n^{ij}(x)=\delta_{ij},\quad w_n^{ij,k}(x) = 0,\quad i,j,k\in \{1,\hdots,d\}\}
\end{align*}
and is therefore measurable with respect to the $\sigma$-algebra $\sigma(h,G,v_n^{ij}, w_n^{ij,k}; i,j,k\in \{1,\hdots,d\})$ even though the partial derivatives  $ \dfrac{\partial u_n^i}{\partial x^l}$ may not be. 
\end{lemma}
\begin{proof} That $\widetilde{Z}_n $ contains $B_n$ follows from the fact that $u_n$ is the identity on $B_n$. As seen from its definition above, the symmetric matrix-valued function $q_n^{\cdot \cdot, k}:=[q_n^{ij,k}]$ can be expressed as  $Du_n \,\text{Hess}(u_n^k)Du_n^T$ where Hess($u_n^k$) is the Hessian matrix $[\frac{\partial^2 u_n^k}{\partial x^i \partial x^j}]$ of second-order partial derivatives. Fix a point $x\in A_n$. As $ h(x) = h\circ u_n(x)$, the full-rank Jacobian matrices $Dh(x) = [\frac{\partial h^i}{\partial x^j}(x)]\; \text{and}\,D(h\circ u_n)(x) =  [\frac{\partial( h^i\circ u_n)(x) }{\partial x^j}] $, acting as matrix coefficients, determine linear injections mapping $\mathbb{R}^d$  onto a common $d$-dimensional tangent space  $T_{h(x)}$ of $h(G)$ at $h(x)$. The  composite linear bijection $\mathbb{R}^d \mapsto T_{h(x)} \mapsto  \mathbb{R}^d $ given by the former followed by the inverse of the latter is represented by the  $d \times d$-matrix $Du_n$, which is therefore nonsingular on $A_n$.
 So the common zero set of the $d$ Hessian matrices Hess($u_n^k$) coincides with that of the $d$  matrix functions $q_n^{\cdot \cdot, k}$. The result then follows from the $\sigma(h,G)$-measurability of $A_n$ and the linear bijection between the families $\{q_n^{ij,k}\}$ and $\{w_n^{ij,k}\}$.
\end{proof}

\begin{lemma}\label{equivalence}
Suppose $h$ and $G$ satisfy the conditions of Proposition \ref{incongruence}. The following three statements are equivalent
\begin{enumerate}[(i)]
\item $h$ possesses a fine local symmetry in $G$.
\item For some $n$, $u_n$ coincides with a nontrivial affine isometry on a finely open set in $A_n \setminus \bar{B}_n$.
\item For some $n$ there is a finely open set $V$ in $Z_n \setminus \bar{B}_n$, where for $d=1 \quad Z_n$  is the set
$\{x \in A_n, \; v_n(x) = 1\}$
and for $d\geq 2\quad Z_n$ is now the set $\{x \in A_n,\; v_n^{ij}(x)=\delta_{ij},\; w_n^{ij,k}(x) = 0,\; i,j,k\in \{1,\hdots,d\}\}$
 introduced in lemma \ref{coincidence}.
\end{enumerate}
Consequently, the fine total asymmetry of $h$ on $G$ is equivalent to the fine closure of $Z_n$ being exactly $\bar{B}_n$ for every $n$.
\end{lemma}

\begin{proof}

Suppose that (i) holds. Then there is a non-trivial isometry $\kappa$ such that $h = h\circ\kappa$ on some finely open set in $S_\kappa$. Choose a sufficiently small finely open bounded set $S\subset S_\kappa$ that does not contain any fixed point of $\kappa$ so that $\kappa(S)$ and $S$ are separated  -- that is, $\inf[|x-y|:x\in S,\;y\in \kappa(S)] >0$ --  and so that $\kappa(S)$ lies in $B_n$ for some $n$. Then $S\subset {A}_n \setminus \bar{B}_n$. Consequently, for $x\in S$,
\begin{align*}
u_n(x) := h_n^{-1}\circ h(x) = h_n^{-1}\circ h\circ \kappa(x) = \kappa(x)
\end{align*}
where the last equality follows from the fact that $h_n^{-1}\circ h(z) = z$ for $z\in \bar{B}_n$. So (i) implies (ii). That (ii) implies (i) follows from the fact that if $u_n$ is an isometry on a finely open set $S\subset{A}_n\setminus \bar{B}_n$, then on $S$, $h\circ u_n = h\circ h_n^{-1}\circ h = h$, and $h$  possesses a local $u_n$-symmetry.
To show that (ii) implies (iii) assume (ii) holds, so that there is a finely open set $S \subset A_n\setminus \bar{B}_n$ on which $u_n$ is an affine isometry.
 If $d=1$ isometries are of the form $\kappa(x) = a \pm x$ and $ S \subset Z_n $ trivially. If $d\geq 2$ isometries are of the form $\kappa(x) = a + Ox$, where $O$ is an orthogonal transformation, and it is clear that $S $ is a subset of the set $\widetilde{Z}_n$, known by Lemma \ref{coincidence} to coincide with $Z_n$.
 Now suppose that (iii) holds and $d=1$. Then $\frac{du_n}{dx}$, being continuous on an open interval $V_0$ in the set $V$ of (iii), is constant at $\pm 1$ on $V_0$, so that $u_n(x)$ is of the form $a \pm x$, implying (ii). For the case where $d\geq 2$ and $V$ is a  finely open set for which (iii) holds, pick a point $x_0\in V$. Then by lemma \ref{GardnerLyons}, there is a finely open neighbourhood $V_0$ of $x_0$ in $V$ such that any point $x$ in $V_0$ is connected to $x_0$ by a polygonal path in $V$. By Lemma \ref{coincidence} $V_0$ lies in $\widetilde{Z}_n$ and $Du_n$ is constant and orthogonal on such sets. Consequently $u_n$ is uniquely defined, given $u_n(x_0)$ and an orthogonal $Du_n(x_0)$, by a double integration along any such rectifiable path. So $u_n(x)$ is of the form $ u_n(x_0) + Du_n(x_0)(x-x_0)$. $Du_n$ being orthogonal, (ii) follows. Thus (i), (ii) and (iii) are all equivalent.
 
 The last assertion of the lemma is just the contrapositive form of the equivalence between (i), (ii) and (iii); that is, the equivalence between the properties: (i') $h$ is finely totally asymmetric, (ii') for every $n,\;B^*_n = A^*_n $; that is, $B_n$ is finely dense in $A_n$ and (iii') for every $n,\;Z_n \setminus B^*_n  $ contains no finely open set, together with the fact that $B^*_n = \bar{B}_n$.

\end{proof}

The proof of Proposition \ref{incongruence} also requires the following lemma, which is a bespoke version of a standard theorem on quadratic variation for continuous semimartingales; our proof of the lemma leans heavily on Protter's proofs of theorems 22 and 32 in Chapter II of \cite{protter2003stochastic}.

\begin{lemma}\label{adapted}
Let $g$ be an $\mathbb{R}^e$-valued Borel function on $\mathbb{R}^d$ that is continuously differentiable on an open set $G$.
  Suppose that $X = (X_t)_{t\geq 0}$ is an $\mathbb{R}^d$-valued Brownian motion with an arbitrary initial distribution $\mu$ on the complete probability space $(\Omega, {\cal X}^\mu, \mathds{P}^\mu)$. Let $\mathbb{G}^\sim$ be the right-continuous filtration $({\cal G}^\sim_{t+})_{t \geq 0}$ where ${\cal G}^\sim_{t}$ is the intersection $\bigcap_{\mu}{\cal G}^{\mu}_{t}$, ${\cal G}^\mu_{t}$ being the $\sigma$-algebra generated by $(g(X_r),\mathds{1}_G(X_r))_{0\leq r\leq t}$ and completed by the $\mathds{P}^\mu$-null sets of ${\cal X}^0_t$. Then, if for $t \geq 0 \quad {\cal G}^\sim_{t,t+}$ denotes the germ field $ \Theta^{-1}_t{\cal G}^\sim_{0+}$, the $e\times e$-matrix-valued process $(Dg(X_t)Dg(X_t)^T\mathds{1}_G(X_t))_{t\geq 0}$ is both $\mathbb{G}^\sim$-progressive and adapted to the process of germ fields $({\cal G}^\sim_{t,t+})_{t\geq 0}$. Furthermore, the matrix-valued function $DgDg^T\mathds{1}_G$ is measurable with respect to the germ enlargement $\widehat{\sigma}(g,G)$ of ${\sigma}(g,G)$.
\end{lemma}

As will be clear from the proof, the matrix-valued process $(Dg(X_t)Dg(X_t)^T\mathds{1}_G(X_t))_{t\geq 0}$ can be viewed as the rate of change of a quadratic-covariation process for $g(X_t)$ that is restricted to random time intervals during which $X_t$ lies in $G$. A simple but notable consequence of this lemma is that for any real-valued function $g$ that is ${\cal C}^1$ on $\mathbb{R}^d$ the process $g \circ X$ possesses a $\mathbb{G}^\sim$-progressive quadratic-variation process with continuously differentiable paths, whether or not it is a semimartingale.
\begin{proof}
We consider the case where $g$ is real-valued; polarization allows us to extend the proof to the multidimensional case. Let $Q_t $ denote the process $|\nabla g(X_t)|^2 \mathds{1}_G(X_t)$. We first establish  the ${\cal G}^\mu_{0+}$-measurability of $Q_0 $ . Notice that 
$\mu$ can always be decomposed as $\mu = \mu(G) \mu' + \mu(G^c)\mu'', \text{ uniquely if } 0 < \mu(G) < 1 $, where $\mu', \; \mu''$ are mutually singular distributions with $\mu'(G) = 1 \text{ and } \mu''(G^c) = 1$. Furthermore $ {\cal G}^\mu_{0+} = {\cal G}^{\mu'}_{0+} \bigwedge {\cal G}^{\mu''}_{0+}$.
Since $ \mathds{P}^{\mu''} (Q_0 = 0) = 1$, $Q_0$ is trivially ${\cal G}^{\mu''}_{0+}$-measurable. Consequently we can restrict our attention to those $\mu$ for which $\mu(G) = 1$. Let $r(x)$ denote the distance  $\inf\{|y-x|: \,y\in G^c\}$. Let $R$ be the  $\mathbb{X}^\mu$-stopping time
\begin{align*}
R = \inf \big\{t\geq 0,\quad \text{either }\, |X_t - X_0| \geq  \frac{r(X_0)}{2} \land 1 \, \text{ or }\, |\nabla g(X_t) - \nabla g(X_0)| \geq 1\big\} \land 1.
\end{align*}
Since we are assuming $\mu(G)=1$ and the paths of $X$ are continuous, $\mathds{P}^\mu(0 < R \leq 1 ) = 1$. Notice that the paths of the stopped process $X^R =(X_{t\land R})_{t \geq 0}$, together with its chords, lie in the convex compact set $\bar{B} (X_0, 1 \land \frac {r(X_0)}{2})$, which is strictly contained in $G$.


Introduce the right-continuous `expanded' filtration due to M. Yor (see \cite{yor1978grossissement} or \cite{protter2003stochastic} p370) that extends the filtration $\mathbb{G}^\mu=({\cal G}^\mu_{t+})_{t\geq 0}$ as follows: $\mathbb{G}^{\mu R} = ({\cal G}^{\mu R}_{t+})_{t \geq 0}$ where ${\cal G}^{\mu R}_{t+} = \{A\in{\cal X}^\mu $ and there exists $ A_t\in {\cal G}^\mu_{t+}:\, A\cap \{R>t\} = A_t\cap \{R>t\}\}$. 
Then $R$ is clearly a $\mathbb{G}^{\mu R}$-stopping time.  As $g \circ X$ has continuous sample paths on $\{ 0\leq t\leq R \}$ and is $\mathbb{G}^{\mu R}$-adapted, the stopped process $g \circ X^R$ is also $\mathbb{G}^{\mu R}$-adapted. We shall first prove that $|\nabla g(X_0)|^2$ is ${\cal G}^{\mu R}_{0+}$-measurable. Let $V_n^R$ be the finite sum
\begin{align*}
V_n^R = \sum_{i=0}^{2^n-1}[g(X_{(i+1)2^{-n}}^R) - g(X_{i2^{-n}}^R)]^2
\end{align*}
Then $V_n^R$ is ${\cal G}_{1+}^{\mu R}$-measurable and can be expressed as
\begin{align*}
V_n^R = \sum_{i=0}^{2^n-1}[(\widehat{\nabla}g_{i,n})^T(X^R_{(i+1)2^{-n}}- X^R_{i2^{-n}})]^2
\end{align*}
where $\widehat{\nabla} g_{i,n}$ denotes the average $\int_0^1 \nabla g((1-\theta)X^R_{i2^{-n}} + \theta X^R_{(i+1)2^{-n}}) d \theta$ of the gradient of $g$ along a chord of $X^R$. Noting that for vectors $a,b$ and $c$
\begin{align*}
(a^Tb)^2 - (c^Tb)^2 = (a+c)^Tbb^T(a-c)
\end{align*}
we see that, almost surely 
\begin{align*}
\left| V_n^R - \right. & \left. \sum_{i=0}^{2^n-1} \right.\left.\big(\nabla g(X^R_{i2^{-n}})\big)^T (X^R_{(i+1)2^{-n}}- X^R_{i2^{-n}}))^2 \right|  \\
&= \left| \sum_{i=0}^{2^n-1}\big(\widehat{\nabla} g_{i,n}+\nabla g(X^R_{i2^{-n}})\big)^T \big(X^R_{(i+1)2^{-n}}- X^R_{i2^{-n}}\big)\times\right.\\
&\left. \hspace{6cm} \big(X^R_{(i+1)2^{-n}} - X^R_{i2^{-n}}\big)^T \big(\widehat{\nabla} g_{i,n} - \nabla g(X^R_{i2^{-n}})\big) \right| \\
&\leq 2\big(|\nabla g(X_0)| +1\big)\left(\sum_{i=0}^{2^n-1}\left| X^R_{(i+1)2^{-n}}- X^R_{i2^{-n}}\right| ^2\right) \left(m_{X_0}^{\nabla g} \circ m^X \left(\frac{1}{2^n}\right)\right)
\end{align*}
where $m_x^{\nabla g}(r)$ is the local modulus of continuity $\sup\{|\nabla g(y) - \nabla g(y')|:|y-y'|\leq r,\; y,y' \in \bar{B}(x,1\land \frac{r(x)}{2})\}$ of $\nabla g$ and $m^X$ is the random modulus of continuity of $X$ on $[0,1]$. The theory of It\^o integration (see Protter \cite{protter2003stochastic} p.77) tells us that the finite series in the first and last expressions in the inequality above both converge in probability as $n\to\infty$. The one on the left converges to $\int_0^R |\nabla g(X_s^R)|^2 ds$ as $n\to\infty$ and the one on the right to $R$. Since $(m_{X_0}^{\nabla g}\circ m^X)(\frac{1}{2^n})\downarrow 0$ a.s.$\mathds{P}^\mu$ as $n\to\infty$, we also have that the integral
\begin{align*}
\int_0^R |\nabla g(X_s)|^2 ds = \int_0^R |\nabla g(X_s^R)|^2 ds = \lim_{n\to\infty}V_n^R,
\end{align*}
the limit being in probability, and is therefore ${\cal G}_{1}^{\mu R}$-measurable for $\mu$ supported on $G$. In the same vein, we can show that $\int_0^{R\land t} |\nabla g(X_s)|^2 ds$ is ${\cal G}_{t}^{\mu R}$-measurable for $t\in[0,1]$. Since, by continuity, 
\begin{align*}
|\nabla g(X_0)|^2 = \lim_{t\downarrow 0}\frac{1}{t}\int_0^{R\land t} |\nabla g(X_s)|^2 ds
\end{align*}
$|\nabla g(X_0)|^2$ is ${\cal G}_{0+}^{\mu R}$-measurable. It remains to prove that ${\cal G}_{0+}^{\mu R} = {\cal G}^\mu_{0+}$. As ${\cal G}^\mu_{0+} \subset {\cal G}_{0+}^{\mu R}$, we shall prove the reverse inclusion. So suppose that $A\in{\cal G}_{0+}^{\mu R}$. Then $A\in{\cal G}_{t+}^{\mu R}$ for $t \geq 0$ and there exists an $A_t\in {\cal G}^\mu_{t+}$ such that $A\cap\{R>t\} = A_t\cap\{R>t\}$. So $\mathds{P}^\mu(A \Delta A_t)\leq 2\mathds{P}^\mu(R\leq t)$. Since $\mathds{P}^\mu(R\leq t)\downarrow 0$ as $t\downarrow 0$, $\limsup_{n} A_{t_n} = A $ a.s., where the limit superior is taken along a sequence $\{t_1, t_2,\hdots\}$ decreasing sufficiently rapidly to zero. Hence 
$A\in \bigcap_{n} {\cal G}^\mu_{t_n+} = {\cal G}^\mu_{0+}$.
The choice of $A$ in ${\cal G}_{0+}^{\mu R}$ is arbitrary; so ${\cal G}_{0+}^{\mu R} \subset {\cal G}^\mu_{0+}$, which establishes the ${\cal G}^\mu_{0+}$-measurability of $|\nabla g(X_0)|^2$. 
Combining this with the earlier result for those $\mu$ for which $\mu(G^c) = 1$, we have also established that $|\nabla g(X_0)|^2 \mathds{1}_G(X_0)$ is measurable with respect to ${\cal G}^{\mu}_{0+}$ for arbitrary $\mu$ and therefore to ${\cal G}^{\sim}_{0+}$. Since $X_t$ can be written $X_0 \circ \Theta_t$, it follows that $(|\nabla g(X_t)|^2 \mathds{1}_G(X_t))$ is adapted to the process of germ fields in the  lemma. As  ${\cal G}^{\sim}_{t+} =\bigcap_\mu\sigma( {\cal G}^{\sim}_{s,s+},\; 0\leq s \leq t, {\cal N}_t^\mu)$, it is also $\mathbb{G}^\sim$-adapted. The measurability of $|\nabla g|^2 \mathds{1}_G$ with respect to $\widehat{\sigma}(g,G)$ follows from the definition of the latter. Polarisation leads to the corresponding results for multidimensional $g$.  This completes the proof.

  \end{proof}

Returning to the proof of  \ref{incongruence}, we now consider the measurability properties of the following processes defined for $t\geq 0$, $n\in\mathbb{N}$ and $i,j,k\in\{1,2,\hdots,d\}$ by
\begin{align*}
F_{n,t} &= \mathds{1}_{A_n}(X_t),\\
U_{n,t}^i &= u_n^i(X_t),\\
V^{ij}_{n,t} &= v^{ij}_n(X_t),\\
W^{ij,k}_{n,t} &= w_n^{ij,k}(X_t), \; \text{ if }\, d \geq 2.
\end{align*}

 Since the functions $\mathds{1}_{A_n}$ and $u_n^i$ are $\sigma(h,G)$-measurable, the processes $F_n$ and $U_{n}^i$ are  $\mathbb{Y}^{\mu}$-progressive for arbitrary $\mu$ and therefore are $\mathbb{Y}^{\sim}$-progressive. Given the continuous differentiability on $A_n$ of the $u_n^i$, lemma \ref{adapted} implies that, for any $\mu$, the processes $V^{ij}_{n}$ are $\mathds{P}^\mu$-indistinguishable from the quadratic-covariation `rate' processes described in that lemma for the  processes $U_{n}^i$; as such, they are $\mathbb{Y}^{\sim}$-progressive as well.  Moreover, for the case where $d\geq 2$, the functions $v_{n}^{ij}$ are also continuously differentiable on $A_n$, and a further application of lemma \ref{adapted}, with $g$ taken as an ordering of
$\{u_n^i,v_n^{ij},\; i,j\in \{1,\hdots,d\},\; i\leq j\}$,
shows that the processes $W^{ij,k}_{n}$ are also $\mathbb{Y}^{\sim}$-progressive. Consequently the functions that are components of $g$, and the $w^{ij,k}_{n}$ where relevant, are all measurable with respect to the germ enlargement $\widehat{\sigma}(h,G)$. The set $Z_n$ defined in lemma \ref{coincidence} is measurable with respect to the $\sigma$-algebra $\sigma(h,G,v_n^{ij}, w_n^{ij,k}; i,j,k\in \{1,\hdots,d\})$ and therefore with respect to the germ enlargement $\widehat{\sigma}(h,G)$. Now by lemma \ref{equivalence} the assumption of fine total asymmetry implies that $\bar{B}_n$ coincides with the fine closure $Z_n^*$, which by lemma \ref{germ h} is also $\widehat{\sigma}(h,G)$-measurable. So $B_n$, as the fine interior of $\bar{B}_n$, is similarly measurable. As a result the map $u_n\mathds{1}_{B_n}$ is $\widehat{\sigma}(h,G)$-measurable for each $n$. But from the construction of $u_n$ this is the local identity id$\mathds{1}_{B_n}$. Since the balls $B_n \in \widehat{\sigma}(h,G)$ and $G = \bigcup_n B_n$, a straightforward piecemeal construction shows that id$\mathds{1}_{G}$ is also $\widehat{\sigma}(h,G)$-measurable. Consequently $\widehat{\sigma}(h,G)$ contains $\widehat{\sigma}(\text{id}\mathds{1}_{G},G)$. But $G$ is finely dense in $\Rd$; so it follows from lemma \ref{dense D} that these two sub-$\sigma$-algebras coincide with ${\cal B}(\Rd)$. Consequently $X_0$ is ${\cal Y}^{\sim}_{0+}$-measurable and it follows from Lemma $\ref{enlargement}$ that  $X$ is adapted to the filtration $\mathbb{Y}^{\sim}$. 
\end{proof}

The following proposition addresses the local problem. It makes precise the intuitive idea that if a Brownian motion starts at a known point and there is a collection of `candidate' $\widehat{\sigma}(h)$-measurable measurement maps that are locally invertible on a sufficiently dense set of points, then the process should be determined at all times by the process of measurements, even if the original measurement $h$ is not finely totally asymmetric. In contrast to the previous proposition, the measurement function $h$ is not assumed to have any properties of smoothness.

A $\sigma$-algebra ${\cal S}$ restricted to a domain $ D \in {\cal B}(\Rd)$ will be said to be  \textbf{locally separable on $D$} if for each $x \in\Rd$ there is a non-empty neighbourhood ball $ B_x$ with centre $x$ such that ${\cal S}$ is separable on $D\cap B_x$; that is, ${\cal S}|_{D\cap  B_x} = {\cal B}(\Rd)|_{D\cap B_x}$.

\begin{proposition}\label{localsuff}
Suppose the $\mathbb{R}^d$-valued Brownian motion $X$ starts at a given point $z$ and $D$ is a $\widehat{\sigma}(h)$-measurable set that is finely dense in $\Rd$. Suppose further that 
the restricted germ enlargement $\widehat{\sigma}(h)|_D$ is locally separable on $D$. Then $X$ is adapted to $\mathbb{Y}^z$.

\end{proposition}

The implication of the condition on the $\widehat{\sigma}(h)$-measurability of $D$ is that `suspect' measurements that lie outside $D$ are recognizable as such and may be ignored with impunity; for instance, if $h(D)\cap h(D^c) = \emptyset$, then $D^c$ might represent, to use D. Rumsfeld's phrase, the set of ``known unknowns".

The proof of Proposition \ref{incongruence} made use of the right continuity of $X$; the following proof of Proposition \ref{localsuff} exploits its local H\"older continuity of order 1/3.

\begin{proof}First, consider some simple consequences of the property of local separability.
\begin{enumerate}[(a)]

\item Lemma \ref{dense D}(i) shows that for any ball $B_x$ the separability of $\widehat{\sigma}(h)$ on $D\cap B_x$ implies that the same property holds on $B_x$. Consequently the unrestricted $\sigma$-algebra $\widehat{\sigma}(h)$ is locally separable on $\Rd$.

\item Writing a neighbourhood ball $B_x$ in the construction above as $B(x, \rho(x))$, we may -- and will -- assume that its radius $\rho(x)$ is strictly positive and is continuous in $x$ . This follows from the fact that the balls $B(x,\rho(x)),\; x\in \Rd$ also form a system of defining  neighbourhood balls if $ \rho(x) = \max_n[(r_n - |x-x_n|)\vee 0]$, where the countable family $\{B(x_n,r_n)\}_n$ is chosen from the given family of defining balls to be a sub-cover of $\Rd$ with finite overlaps at each point.

\end{enumerate}

 Given the continuity of the paths of $X$ and the right-continuity of $\mathbb{X}^z$ and $\mathbb{Y}^z$, to prove the proposition we need only consider the measurability of $X_t$ on the dyadic rationals $\mathbb{D}_+ =\cup_n\mathbb{D}_n$, where $\mathbb{D}_n = \{t \geq 0:\;2^nt \in \mathbb{N}\}$.
That is, we have to establish that for $t \in \mathbb{D}_+, \;{\cal Y}^z_{t+} = {\cal X}^z_{t}$. Given that $X$ is strong Markov, it enough to show that $X_t$ is  ${\cal Y}^z_{t+}$-measurable. This is equivalent to the statement
\begin{equation}\label{eq3}
\mathds{P}^z(  f_i(X_{t}) = \mathds{E}^z[f_i(X_{t})| {\cal Y}^z_{t+}]) = 1, \qquad t\in \mathbb{D}_+
\end{equation} 
for each member $f_i$ of a countable collection of $d$ strictly positive bounded Borel functions on $\Rd$ that generate ${\cal B}(\Rd)$. To prove this we introduce a nondecreasing sequence of $\mathbb{X}^z$-adapted stopping times $(T_n)_n$ that is localizing in the sense that for every $t \in \mathbb{D}_+ \;\;\mathds{P}^z(t < \lim_n T_n ) =1 $, for which the following is satisfied: for each $f_i,\;n$ and $t\in \mathbb{D}_n$
\begin{equation}\label{eq4}
\mathds{P}^z(f_i(X_{t})\mathds{1}_{\{t<T_n\}} =  \mathds{E}^z[f_i(X_{t})\mathds{1}_{\{t<T_n\}}| {\cal Y}^z_{t+},t<T_n]) = 1.
\end{equation} 
Since $f_i(X_{t})\mathds{1}_{\{t<T_n\}}$ is only zero on $\{t\geq T_n\}$, this equation can be read as:
\begin{equation}
 For\; each\; n \; and \; t\in \mathbb{D}_n, \;\; X_t,\; when\; restricted\; to\; \{\omega:t<T_n(\omega)\}, \;\; is\;\; {\cal Y}^z_{t+}|_{\{t<T_n\}}-measurable.\nonumber 
 \end{equation}
Since we are assuming that $T_n \nearrow \infty \;a.s. \mathds{P}^z$ with $n$, we have that ${\cal Y}^z_{t+} = \sigma(\cup_n{\cal Y}^z_{t+}|_{\{t<T_n\}})$. Consequently equation (\ref{eq3}) would follow from (\ref{eq4}). 
So it remains to construct the $T_n$ so that (\ref{eq4}) is valid.

First introduce the dyadic-valued stopping times $R_n' \text{ and }S'_n  $ where $R_n' = \min\{t \in \mathbb{D}_n:\rho(X_t)\leq 2^{-n/3} \}$ and $S'_n = \min\{t \in \mathbb{D}_n:M^n_t\geq 2^{-n/3} \}$, $M^n_t$ being the modulus of consecutive increments before $t: \;\max\{|X_s-X_{(s-2^{-n})\vee 0}|:s\in \mathbb{D}_n\cap[0,t]\}$.
These are adapted to the sub-filtration $ \mathbb{X}^{z,n} = (\mathcal{X}^z_t)_{t\in \mathbb{D}_n }$. Then set $T'_n = R'_n\wedge S'_n$, noting that $T'_n > 0 $ if and only if $\rho(z)>2^{-n/3}$.

It can be checked that $(X_t)_{t \in \mathbb{D}_n\cap [0, T'_n)}$ is a discrete-parameter sub-Markov process, restricted as it is to the $\mathbb{X}^{z,n}$ stochastic interval $\Lambda_n = \{(t,\omega): t< T'_n(\omega) \} \subset \mathbb{D}_n \times \Omega$. This can be made Markov, in the manner of Doob, by expanding the state space of $X$ to the Borel space $(\Rd_\partial =  \Rd \cup \{\partial\},{\cal B}(\Rd_\partial)$ and regarding the additional point $\partial$ as an isolated `trap' state to which the process jumps at $T'_n$. Let $X^{n\partial}_t$ be the process that is $X_t$ on $\{t< T'_n\}$ and $\partial$ on $\{t \geq T'_n\}$. Then $(X^{n\partial}_t)_{t \in \mathbb{D}_n}$ is a Markov process stopped at the point $\partial$. Note that $X^{n\partial} $ coincides with $X$ on $\Lambda_n$ and only takes the value $\partial$ on the complementary set $\mathbb{D}_n \times \Omega \setminus \Lambda_n$. Note the equivalences $\{t< T'_n\} = \{X^{n\partial}_t \in \Rd\} = \{X^{n\partial}_t = X_t\}$.
Now consider the hybrid sub-$\sigma$-algebra in ${\cal X}^z$:
\begin{equation}
{\cal Y}^{z,n}_t = \sigma(X^{-1}_s \widehat{\sigma}(h)|_{\{s< T'_n\}}: s\in \mathbb{D}_n\cap(0,t],\;{\cal N}^z_t ), \quad t \in \mathbb{D}_n.
\end{equation}
This is generated by the `germ' observations of $X_s$ at those dyadic times of order $n$ that are before $T'_n$ and not beyond $t$. Since for all $s  \in \mathbb{D}_n,\;X^{-1}_s \Rd \cap{\{s< T'_n\}} =\{s< T'_n\},\;T'_n $ is a $({\cal Y}^{z,n}_t)_{t  \in \mathbb{D}_n}$-stopping time. 
 Let $\widehat{\sigma}_\partial(h) = \sigma(\widehat{\sigma}(h),\{\partial\})$ be the natural extension of $\widehat{\sigma}(h)$ to a sub-$\sigma$-algebra of ${\cal B}(\Rd_\partial)$ that includes $\{\partial\}$.
  Then, in terms of $X^{n\partial}, \;{\cal Y}^{z,n}_t$ can be written as $\sigma(X^{n\partial,-1}_s \widehat{\sigma}_\partial(h): s\in \mathbb{D}_n\cap(0,t],\;{\cal N}^z_t )$; that is, the corresponding completed filtration is induced by the map $X^{n\partial}$ on $\mathbb{D}_n \times \Omega$ into the `observation' space $(\Rd_\partial,\widehat{\sigma}_\partial(h))$.

In the next two paragraphs, the dyadic increment $2^{-n}$ has been shortened to $\Delta $ for clarity, the dependence on $n$ being understood.
The next step is to establish that for $t \in \mathbb{D}_n,\; X_t|_{\{t< T'_n\}}$ is $ {\cal Y}^{z,n}_t $-measurable or, equivalently, that $X^{n\partial}_t$ is ${\cal Y}^{z,n}_t$-measurable. We shall prove this by induction along $\mathbb{D}_n$. 
First modify the strictly positive bounded functions $f_i $ that generate ${\cal B}(\Rd)$, so that they generate ${\cal B}(\Rd_\partial)$ as well, by setting $f_i(\partial) = 0$ for all $i$. Assume for some $t>0$ that $X^{n\partial}_s$ is ${\cal Y}^{z,n}_s$-measurable for $s\in \mathbb{D}_n\cap(0,t-\Delta]$. Then for each $f_i $ and $n$, up to $\mathds{P}^z$ null sets,

\begin{eqnarray}
\lefteqn{\mathds{E}^z[f_i(X_{t})\mathds{1}_{\{t<T'_n\}}| {\cal Y}_t^{z,n}]}\nonumber\\
  & = & \mathds{E}^z[f_i(X^{n\partial}_{t}) |{\cal Y}^{z,n}_t]\nonumber\\
  & = & \mathds{E}^z[f_i(X^{n\partial}_{t})|{\cal Y}^{z,n}_{t-\Delta},\;X^{n\partial,-1}_{t}\widehat{\sigma}_\partial(h) ]\nonumber\\
  & = & \mathds{E}^{X^{n\partial}_{t-\Delta}}[f_i(X^{n\partial}_{t})\mathds{1}_{\Rd}(X^{n\partial}_{t})|\;X^{n\partial,-1}_{t}\widehat{\sigma}_\partial(h) ]\nonumber\\
  & = & \mathds{E}^{X^{n\partial}_{t-\Delta}}[f_i(X^{n\partial}_{t})\mathds{1}_{\Rd}(X^{n\partial}_{t-\Delta}) \mathds{1}_{C}(X^{n\partial}_{t-\Delta}, X^{n\partial}_{t})|\;X^{n\partial,-1}_{t}\widehat{\sigma}_\partial(h) ]\nonumber\\  
  & = &  \mathds{1}_{\Rd}(X^{n\partial}_{t-\Delta}) \mathds{E}^{X^{n\partial}_{t-\Delta}}[ f_i(X^{n\partial}_{t}) \mathds{1}_{C}(X^{n\partial}_{t-\Delta}, X^{n\partial}_{t}) ] \nonumber\\
  & = &  f_i(X^{n\partial}_{t})\mathds{1}_{\Rd}(X^{n\partial}_{t})\nonumber\\ 
  & = & f_i(X_t) \mathds{1}_{\{t<T'_n\}}\qquad a.s.\;\mathds{P}^z
\end{eqnarray} 
The first three equations follow from the relations $\{ X^{n\partial}_t \in \Rd \} = \{ X^{n\partial}_t = X_t \} = \{t<T'_n\} \subset \{t-\Delta<T'_n\} = \{ X^{n\partial}_{t-\Delta} \in \Rd \} $, expressed as indicator functions, together with the inductive assumption and the Markov property of $X^{n\partial}$.
 In the fourth and fifth equations, the set $\{ X^{n\partial}_t \in \Rd \}$ has been decomposed as $\{X^{n\partial}_{t-\Delta} \in \Rd\}\cap \{(X_{t-\Delta}, X_t) \in C\}$ where  $C =\{ (x,y)\in \Rd \times \Rd :\;\rho(x)\wedge\rho(y)> \Delta^{\frac{1}{3}}, \;|y - x|< \Delta^{\frac{1}{3}} \}$, which follows from the definition of $T'_n$ as $R'_n\wedge S'_n$. However, by construction, $\widehat{\sigma}(h)$ is separable on each of the Borel $x$-sections of $C$. Consequently for $X^{n\partial}_{t-\Delta} \in \Rd$, the distribution of $X^{n\partial}_{t}$, conditioned on $X^{n\partial}_{t-\Delta}$, is a Dirac measure, which justifies the sixth equation. 
 The remaining equation is a reversion to the original notation. 
 
 Since for $t=0, \;\;\mathds{E}^z[f_i(X_{0})\mathds{1}_{\{0<T'_n\}}| {\cal Y}_0^{z,n}] = f_i(z)\mathds{1}_{\{0<T'_n\}}$, we can apply this result inductively to establish that for $t \in \mathbb{D}_n,\; X_t|_{\{t<T'_n\}}$ is ${\cal Y}^{z,n}_t$-measurable.

 Since ${\cal Y}^{z,n}_t|_{\{t<T'_n\}} \subset {\cal Y}^{z}_{t+}|_{\{t<T'_n\}}$, it follows that 
 \begin{equation}
 \mathds{P}^z(f_i(X_{t})\mathds{1}_{\{t<T'_n\}} =  \mathds{E}^z[f_i(X_{t})\mathds{1}_{\{t<T'_n\}}| {\cal Y}^z_{t+},t<T'_n]) = 1;
 \end{equation}

\noindent that is, (4) holds if the stopping time  $T_n$ is taken to be $T'_n$. But it also holds  for any $[0,\infty]$-valued stopping time $T$ provided $T \leq T'_n $ a.s.$\mathds{P}^z$. We can construct a localising sequence $(T_n)_n$ of such stopping times as follows. Set $T_n = R_n \wedge S_n$ where 
 $R_n =\inf \{ s \geq 0:\rho(X_s) \leq 2^{-n/3} \}$ and $S_n = \inf\{ t:\sup_{\delta \leq 2^{-n}}\frac{M_t(\delta)}{\delta^{1/3}} \leq 1 \}$. Here $M_t(\delta)$ is the random modulus of continuity $\sup\{|X_s - X_{s'}|:\; s,s' \in [0,t],\;  |s-s'|\leq \delta\}$. Levy's theorem on the modulus of continuity of Brownian motion (\cite{revuz2013continuous} p.29) implies that for each $t$ the H\"older-weighted modulus  $\frac{M_t(\delta)}{\delta^{1/3}}$ satisfies $\mathds{P}^z(\lim_{\delta \downarrow 0}\frac{M_t(\delta)}{\delta^{1/3}} = 0) = 1$. Clearly $R_n \leq R'_n$ and $S_n \leq \bigwedge_{m\geq n} S'_m$ and both $(R_n)_n$ and $(S_n)_n$ are non-decreasing localising sequences. So $(T_n)_n$ has the same property and (3) then follows from (4). This completes the proof.
\end{proof}

\section{Applications}
The following justifies the assertion made in the introduction about a real-valued Brownian motion with a real-analytic measurement function. Note that $h$ is said to be \textit{asymmetric} if there are no nontrivial affine isometries $\kappa$ for which $h \equiv h \circ \kappa $ and that the only affine isometries on the real line are  translations and reflections about arbitrary points.
\begin{proposition}\label{realanalytic}

 Suppose $X$ is an $\Rd$-valued process locally equivalent in law to a Brownian motion with a strictly positive initial density. Suppose $h:\Rd \mapsto \mathbb{R}^e,\;e\geq d$ is real-analytic and $Dh(x)$ is of full rank $d$ at some point $x_0$. Then the following statements are equivalent:
\begin{enumerate}[(i)]
\item $h$ is asymmetric,
\item $h$ is finely totally asymmetric,
\item $X$ is $\mathbb{Y}^\sim$-adapted.
\end{enumerate}
If $X$ is of the form $X_t= \Psi \circ X'_t, \; t \geq 0,$ where  $ X'$ is locally equivalent in law to an $\Rd$-valued Brownian motion with a positive initial density and $\Psi $ is a real-analytic diffeomorphism on $\Rd$, then these statements are still equivalent if $h$ is replaced by $h \circ \Psi $.
\end{proposition}

The condition that the initial distribution $\mu$ possesses a density could be weakened to ``$\mu$ charges all finely open Borel sets''; that is, $\mu(G) > 0$ for all finely open Borel $G$.  

\begin{proof} The rank condition implies that the positive-semidefinite symmetric $d \times d$ matrix $J = Dh^tDh$ is such that $\det{J(x_0)}>0$ for some $x_0$. So the function $\det{J}$, being real-analytic, is non-zero everywhere except on its closed zero set $\{ x : \det{J(x)} = 0\} $, which is necessarily $\lambda_d $-null. If we take $G$ to be its complement $\{ x : \det{J(x)} > 0\} $, Corollary \ref{simplification} tells us that the last two statements are equivalent. Clearly $(ii)$ implies $ (i)$; so to complete the proof it is sufficient to show that $(i)$ implies $ (ii)$. Suppose, contrapositively, that $h$ has a fine local symmetry; that is, for some non-trivial affine isometry $\kappa$ the symmetry set $S_\kappa = \{x: h(x)= h\circ \kappa (x) \}$ has a nonempty fine interior $S^{fo}_\kappa$. Pick an interior point $x$ in this set. By Lemma \ref{GardnerLyons} within $S^{fo}_\kappa$ there is a finely connected finely open neighbourhood $V$ of $x$ such that if $y \in V$ and $H_{x,y}$ is the bisecting hyperplane  $\{ z: |z-x| = |z-y| \}$ of the line segment $[x,y]$ there is a set $A \subset H\cap V$ such that the line segments $[x,z]$ lie in $V$ for all $z \in A$; furthermore  $\lambda_{d-1}(A) >0$. Consequently we can choose $d$ linearly independent line segments $[x,z_i]$ lying in $V$ along which the directional derivatives of the real-valued functions $g^j = h^j-h^j \circ \kappa$ are all zero at $x$. So their partial derivatives $\partial g^j/\partial x^i$ are also zero at $x$. 
So proceeding inductively we can show that the functions $g^j$ and all their partial derivatives of all orders are zero at $x$. 
 But the $g^j$ are real-analytic on $\mathbb{R}^d$ and so must be zero everywhere; that is, $S_\kappa = \mathbb{R}^d$. Thus $h$ is symmetric for the isometry $\kappa$. The last statement on the change of coordinates is self-evident.
\end{proof}

The following corollary to Proposition \ref{localsuff} addresses the local problem for a  class of one-dimensional diffusions observed through piecewise-monotone measurement maps. It can be regarded as a `noiseless' supplement to the more specific asymptotic result in \cite{fleming1989piecewise} for unimodal measurement maps with vanishingly small measurement noise. In \cite{fleming1989piecewise} Fleming and Pardoux essentially solve the global problem for real-valued diffusions with $C^2$ coefficients and a $C^2$ unimodal measurement function as a simple by-product of their detailed analysis of the rate of convergence of a family of approximate filters as the measurement noise vanishes. They also solve, in essence, the partial local problem of continuously tracking a diffusion from a known starting point up to the first passage time of the critical point of $h$. The following proposition, which generalizes this last result, solves the full local problem of tracking a real-valued diffusion through all the critical points of a piecewise monotone measurement function $h$.

To state the result it is necessary to introduce some further notation, which we adapt from   \cite{fleming1989piecewise}. Suppose $h$ is a real-valued  $C^1(\mathbb{R}^1)$ function that is strictly monotone between isolated critical points (that is, points $c$ where $h'(c)=0$ and $h(c)$ is a local extremum). For such a function, for each critical point $c$ an open interval $(l,r)$  can be chosen containing $c$ but no other critical points, such that $h|_{[l,r]}$ is unimodal with $h(l) = h(r)$. Let $h_l$ and $h_r$ denote the corresponding strictly monotone restrictions $h|_{[l,c]}$ and $h|_{[c,r]}$, which share a common range $ h ([l,r])$. Finally let $x^-_c= h_l^{-1}\circ h \;\text{and}\; x^+_c = h_r^{-1}\circ h$ be the `left' and `right' functions mapping $[l,r]$ onto $[l,c]$ and $[c,r]$. Note that
, for $x \neq c, \;\{x^-_c(x),x^+_c(x)\}$ are the two-point level sets of $h|_{[l,r]}$ consisting of $x$ and its `alias'.

\begin{proposition}\label{coroldiffusion} 

Suppose $X  = (X_t)_{0\leq t < \infty}$ is a real-valued diffusion satisfying the It\^o equation
\begin{equation*}
dX_t = f(X_t)dt + g(X_t)dW_t, \quad X_0 = z, \text{ a constant}
\end{equation*}
where  (i) $f$ is a Borel function satisfying  a linear growth condition: $|f(x)|\leq K(1 +|x|) $ for some $K$, (ii) $g$ is a positive real-valued  $C^1(\mathbb{R}^1)$ function such that $g,1/g$ and $g'$ are bounded, and (iii) $W$ is a Brownian motion adapted to $\mathbb{X}^z$.  Suppose h is a real-valued $C^1(\mathbb{R}^1)$ function with isolated critical points. Suppose for each critical point $c$ the associated function $\gamma_c =  (gh')^2\circ x^-_c -(gh')^2\circ x^+_c$, defined on a corresponding neighbourhood interval $[l_c,r_c]$, is non-zero on a dense subset. Then $X$ is adapted to $ \mathbb{Y}^z$, the right-continuous filtration generated by the process $Y_t = h(X_t),\;0 \leq t < \infty$ and locally augmented with the null sets of $\mathbb{X}^z$.

\end{proposition}

The condition on the function $\gamma_c$ is essentially a weakening of a ``detectability condition'' $(2.1)$ to be found in \cite{fleming1989piecewise}. It would be straightforward to extend the proposition to include discontinuities in $h$ and $h'$, but only at the expense of an even more cumbersome notation.

\begin{proof} 
First consider the special case where $g \equiv 1$; that is, $X$ is a drifting Brownian motion, defined on the canonical space of its sample paths $(C[0,\infty),\mathcal{X}^0_\infty) $, that satisfies the It\^o equation $dX_t = f(X_t)dt + dW_t, \; X_0 = z,$.
 Girsanov's theorem implies that, since $f$ has linear growth, the probability law defined by the weak solution of this equation -- that is, the law $\mathds{Q}^z$ of $X$ on $\mathcal{X}^0_\infty$ -- exists and is unique. Furthermore, for each $z$ and $t>0$, $\mathds{Q}^z$ is equivalent to the law $\mathds{P}^z$ of Brownian motion starting at $z$  when both are restricted to $\mathcal{X}^0_{t+}$; that is, the two restricted laws share the same null sets and the completion $\mathcal{X}^z_{t+}$ is uniquely defined. Consequently we can apply Proposition \ref{localsuff}, provided the conditions on $h$ are met. Notice that by Lemma 3.4  $(h'(X_t))^2)_{0 \leq t<\infty} $ is $\mathbb{Y}^\sim$-adapted, as it is the rate of change of the quadratic variation process of $Y$. Consequently the function $x \to h'(x)^2$ is  $\widehat{\sigma}(h)$-measurable.
We shall prove that for each real point $u$ the $\widehat{\sigma}(h)$-measurable map $ x \to \bar{h}(x) = ( h(x),h'(x)^2) $ is invertible on a dense subset of an interval neighbourhood $[l_u,r_u]$  of $u$. (Remember that in one dimension the Euclidean and fine topologies are the same). For points $u$ in the open interval $(c_l,c_r)$ between two neighbouring critical points (and the points $\pm \infty$ where appropriate) $h$ is strictly monotone and an inverse map $\bar{h}^{-1}:\bar{h}((c_l,c_r)) \mapsto (c_l,c_r)$ can be defined as $\bar{h}^{-1}(y,v) = h^{-1}(y)$ with values in that interval.
 So consider neighbourhoods for the countable family of critical points ${ C}$. From the conditions of the proposition we can choose for each $c \in C$ a neighbourhood $(l_c,r_c)$ with the properties above that is sufficiently small to be disjoint from the corresponding neighbourhoods of adjacent critical points, so that the set $ D_c = \{x:\gamma_c(x) \neq 0\}\cap (l_c,r_c)$ is dense in $ (l_c,r_c)$. Then for each $c$ the inverse map $\bar{h}_c^{-1}$ is well defined on $\bar{h}(D_c)$ and can be expressed as 
\begin{align*}
\bar{h}_c^{-1}(y,v)   &= h_{l_c}^{-1}(y)\quad    \text{if}\quad v = ( h' \circ h_{l_c}^{-1}(y))^2,\\
                        &= h_{r_c}^{-1}(y)\quad    \text{if}\quad v = ( h' \circ h_{r_c}^{-1}(y))^2.
\end{align*} 

and $\bar h$ is invertible on $D_c$. Consequently the $\widehat{\sigma}(h)$-measurable map $\bar h$ is invertable on a neighbourhood of every real point and $\widehat{\sigma}(h)$ is locally separable. The conditions of Proposition \ref{localsuff} are thus fulfilled and, for the case where $g \equiv 1$, $X$ is adapted to $\mathbb{Y}^z$.

For the general case we begin by transforming $X$ into a process with a law that is locally equivalent to that of Brownian motion.
Introduce the real-valued map $\Phi(x) = \int_{[0,x]}1/g(x')dx'$. Let $\mathring{X}$ be the transformed process $\mathring{X}_t = \Phi(X_t),\; 0 \leq t < \infty$.
The conditions on $g$ show that both $\Phi$ and its inverse are strictly increasing $C^2$ Lipschitz functions. So it is sufficient to show that $\mathring{X}$ is adapted to the filtration $\mathbb{Y}^z$.
Note the observation process can be written as $Y_t = \mathring{h}(\mathring{X}_t), 0\leq t < \infty$ where $\mathring{h} = h \circ \Phi^{-1}$.

 But it follows from It\^o's rule that $\mathring{X}$ satisfies the equation
  \begin{equation*}
  d\mathring{X}_t = b(\mathring{X}_t)dt + dW_t, \; \mathring{X}_0 = \Phi(z),\; 0 \leq t < \infty
  \end{equation*}
 
 where  $b = [f/g  - g'/(2g^2)] \circ \Phi^{-1}$, which is also a Borel function satisfying a linear growth condition. Consequently the formulation is reduced to the special case already considered in the first paragraph, with the measurement map $h$ replaced by $h\circ \Phi^{-1}$ and its first derivative $h'$ by $(gh')\circ \Phi^{-1}$. This completes the proof.
\end{proof}

\bigskip

\section{Appendix}

We include here the justification of Example \ref{example} which contains an $\Rd$-valued Brownian motion tracked from a real-valued observation process.
In the example, we choose $X = (X_t)_{t\geq 0}$ to be an $\Rd$-valued Brownian motion on the usual canonical space and $h$ to be a real-valued function on $\Rd$ of the form

$$ h(x) = e^{a_1x_1} + \ldots + e^{a_dx_d},\ \ x=(x_1,\ldots, x_d)\in\Rd,
$$

\noindent where the constants $a_1, a_2,\ldots, a_d$ are arbitrarily chosen coefficients subject only to the condition that the d-vector $(a_1^{-2},\ldots,a_d^{-2})$ does not lie on any of the $3^d/2$ hyperplanes given by $ \sum_1^d u_jx_j = 0, $ where the $u_j \in \{ -1, 0, 1\}.$ 
For such an $h$, suppose $Y$ is the continuous real-valued process $Y_t = h(X_t)_{ t\geq 0}$. Remember that ${\cal Y}^\mu_{t+}$ denotes the $\sigma$-algebra $\bigcap_{\; s \geq t} \sigma^\mu (Y_r:\; 0 \leq r \leq s, {\cal N}^\mu_s)$ and $\mathbb{Y}^{\sim}$ denotes the right-continuous filtration $(\bigcap_{\mu} {\cal Y}^\mu_{t+})_{t \geq 0}$.
 
 \begin{proposition}
 $ X$ is adapted to the filtration $\mathbb{Y}^{\sim}$ generated by $Y$.
 \end{proposition}
 
 \begin{proof} Let $h_1 = h$. First introduce the function $h_2 = \sum_1^d\langle \partial_jh,\partial_jh \rangle = a^2_1e^{2a_1x_1} + \ldots + a^2_de^{2a_dx_d}$, and for $n = 3, 4, \dots $, the functions
$$h_n = \frac{1}{n-1}\sum_j a_j^{2n-2}\langle \partial_jh,\partial_jh_{n-1} \rangle = a^{2n-2}_1e^{na_1x_1} + \ldots + a^{2n-2}_de^{na_dx_d}.$$

Let $Y_t^{(n)}$ denote the  continuous process $h_n(X_t)$. Since $h$ is smooth, It\^o calculus shows that  $\int_0^t Y_s^{(n)}ds$ is indistinguishable from the continuously differentiable quadratic covariation process $\frac{1}{n-1}\langle Y,Y^{(n-1)}\rangle_t$ and that each $Y^{(n)}$ is adapted to the germ field filtration of $Y$ for every $\mu$. Consequently each $h_n$ is measurable with respect to the `germ field projection' $\widehat{\sigma}(h):= \{ A \in {\cal B}(\Rd):\; \{X_0 \in A\} \in \bigcap_\mu {\cal Y}^\mu_{0+}\}$.
It is notationally convenient to change to positive coordinates. For $j\in \{1,\ldots,d\}$ let $f_j(x_j)= a_j^2e^{a_jx_j}$; the vector function $f$ is a diffeomorphic isomorphism of $\Rd$ onto $(0,\infty)^d $. For $w\in (0,\infty)^d$ let  $\tilde{h}_n(w) = a_1^{-2}w_1^n + \ldots + a_d^{-2}w_d^n$ for $n = 1,2,\ldots.$  Note that $\;h_n = \tilde{h}_n\circ f$ and that each $\tilde{h}_n$ is measurable with respect to ${\cal H}:= f \widehat{\sigma}(h)$, the germ field projection `pushed forward' into ${\cal B}((0,\infty)^d)$. 
We aim to prove that $\widehat{\sigma}(h)= {\cal B}(\Rd)$ or, equivalently, that ${\cal H} = {\cal B}((0,\infty)^d)$. This is the same as establishing that the identity on $(0,\infty)^d$ is ${\cal H}$-measurable. From the $\widehat{\sigma}(h)$-measurability of the functions $h_n$ it is clear that the functions $\tilde{h}_n$ are also ${\cal H}$-measurable.

 Let $\bar{w} = \{\bar{w_1}, \ldots,\bar{w_d}\}$ denote the following reordering of the components of $w \in (0,\infty)^d$:

$\bar{w_1} := \max\{w_1,\ldots,w_d\} \text{ and for } k = 2,\ldots,d,\; \bar{w_k} := \max\{ w_j\mathds{1}_{\{j:w_j < \bar{w}_{k-1}\}}\}$.

\noindent That is, $\bar{w}$ lists in decreasing order the $m\;(\leq d)$ distinct values of the positive $w_j$, followed by $d-m$ zeroes. Let $[c_{jk}(w)]$ denote the $d\times d$-matrix-valued function of zeroes and ones given by $c_{jk}(w) = \mathds{1}_{\{w_j = \bar{w}_k\}}$. Note that each `row' vector $c_{j\star}$ contains a single component of value `one' and that the `column' vectors $c_{\star k}$ are mutually orthogonal, being indicator functions of the disjoint `level' sets in $\{1,\ldots,d\}$ of $w$ (a vector of all zeroes being the indicator function of the empty set).

Now consider the $(0,\infty)^d$-valued map $W(w):= \sum_k \bar{w_k}(w)c_{\star k}(w)$. Closer inspection of this equation will reveal that $W$ coincides with the identity map on $(0,\infty)^d$. It remains to prove that $W$ is ${\cal H}$-measurable or, equivalently, that the  functions $\bar{w_k}(w) \text{ and the vectors }c_{\star k}(w)$ are ${\cal H}$-measurable. Since
$ \max\{w_1,\ldots,w_d\} = \lim_n (\tilde{h}_{n+1} / \tilde{h}_n) $, 
$\bar{w_1}$ is ${\cal H}$-measurable, as is the related partial sum function
$b_1 := \sum_{j:w_j = \bar{w_1}}a^{-2}_j =\lim_n (\tilde{h}_n/ \bar{w_1}^n) $.
Let $\bar{b}:= \sum_j a_j^{-2}$. For $k = 2,3,\ldots,d$, the ${\cal H}$-measurability of the $\bar{w_k}$ and the associated partial sums
$b_k$ follows from that of the right-hand members of the following equations, when successively applied:
\begin{align*}
\bar{w_k} &:= \max\{ w_j\mathds{1}_{w_j < \bar{w}_{k-1}}\} =  \mathds{1}_{\{0<b_{k-1}<\bar{b}\}} \lim_n[ (\tilde{h}_{n+1} - \sum_{j=1}^{k-1} b_j \bar{w_j}^{n+1})/( \tilde{h}_n -\sum_{j=1}^{k-1} b_j \bar{w_j}^{n})],\\ 
b_{k} & := \sum_ja^{-2}_j \mathds{1}_{\{j:\;w_j = \bar{w_k}\}}= \mathds{1}_{\{0<b_{k-1}<\bar{b}\}} \lim_n [( \tilde{h}_n -\sum_{j=1}^{k-1} b_j \bar{w_j}^{n})/ (\bar{w_k}^n)].
\end{align*}

Introduce the linear functional $\beta(z) := \sum_j a_j^{-2}z_j$ on $[0,\infty)^d$ and let $\beta_Q$ denote its restriction to the subset of `indicator' vectors $Q:= \{0,1\}^d.$
The $3^d/2$ constraints on the $a_j$ ensure that $ \quad\beta_Q(c) \neq \beta_Q(c')$ for different $c,\;c' \in Q$; that is, $\beta_Q$ maps the $2^d$ elements of $Q$ onto $2^d$ distinct  points in $[0,\infty)$. Consequently $\beta_Q^{-1}$ is a point-to-point map and $\beta_Q^{-1} \circ \beta_Q$ is the identity on $Q$. Since $b_k = \beta_Q (c_{\star k})$, it follows that the vector functions $\beta_Q^{-1} \circ b_k(w) = \beta_Q^{-1} \circ \beta_Q (c_{\star k}) = c_{\star k}$ are also ${\cal H}$-measurable. This establishes the ${\cal H}$-measurability of the identity map $W(w) = w =  \sum_k \bar{w_k}(w)c_{\star k}(w)$. Consequently $\widehat{\sigma}(h) = {\cal B}(\Rd)$ and Lemma \ref{enlargement} $ X$ is adapted to $\mathbb{Y}^\sim$. \end{proof}

\bibliographystyle{alpha}
\bibliography{biblio}
\end{document}